\documentclass{jktrX}

\def\draftnote{ }
\def\today{}
\arraycolsep=\labelsep

\def\marusen{\unitlength.1em
  \begin{minipage}{10\unitlength}
    \begin{picture}(10,10)
      \put(5,6){\circle{6}} 
      \qbezier(3,8)(7,4)(7,4)
    \end{picture}
  \end{minipage}
}
\def\marubackslash{{\small \marusen}}

\def\hsymb#1{\mbox{\strut\rlap{\smash{\large$#1$}}\quad}} 
\def\Hsymb#1{\mbox{\strut\rlap{\smash{\Large$#1$}}\quad}} 

\usepackage{amssymb}
\usepackage{amsmath}
\usepackage{latexsym}
\usepackage{graphicx}

\renewcommand\hat{\widehat}
\renewcommand\tilde{\widetilde}

\newcommand{\bysame}{\mbox{\rule{3em}{.4pt}}\,}

\def\qed{\hfill \fbox{}}

\newtheorem{thm}{Theorem}[section]

\newtheorem{lemm}[thm]{Lemma}

\newtheorem{prop}[thm]{Proposition}
\newtheorem{cor}[thm]{Corollary}

 \newtheorem{dfn}[thm]{Definition}
\newtheorem{rem}[thm]{Remark}
 \newtheorem{ex}[thm]{Example}

\newcommand\ZZ{{\mathbb Z}}

\newcommand\CC{{\mathbb C}}
\newcommand\QQ{{\mathbb Q}}

\newcommand\nana{
\begin{picture}(12,12)
\multiput(1,0)(5,4){3}{\circle*{1}}
\end{picture}
}

\newcommand\mtx[4]{
\left[\begin{array}{cc}
#1&#2\\
#3&#4
\end{array}\right]}

\newcommand\susum[2]{\textstyle\sum\limits_{#1}^{#2}}

\newcommand\mtxc[4]{
\left[\begin{array}{cc}
#1&#2\\
#3&#4
\end{array}\right]}

\begin{document}

\pagestyle{myheadings}
 \markboth{
Twisted Alexander polynomials  
associated to metacyclic representations
}
{
 \ \ M. Hirasawa \& K. Murasugi}

\title{Twisted Alexander polynomials of $2$-bridge knots 
associated to metacyclic representations}
\author{Mikami Hirasawa}
\address{Department of Mathematics,
Nagoya Institute of Technology\\
Nagoya Aichi 466-8555 Japan\\
{\it E-mail: hirasawa.mikami@nitech.ac.jp}
}
\author{Kunio Murasugi}
\address{Department of Mathematics,
University of Toronto\\
Toronto, ON M5S2E4 Canada\\
{\it E-mail: murasugi@math.toronto.edu}
}

\maketitle


\begin{abstract}
Let $p=2n+1$ be a prime and $D_p$ a dihedral group of order $2p$. 
Let $\widehat{\rho} : G(K) \rightarrow D_p \rightarrow 
GL(p,\ZZ)$ be
a non-abelian representation of the knot group $G(K)$ of a knot $K$ in
3-sphere. 
Let $\widetilde{\Delta}_{\widehat{\rho},K} (t)$ be the twisted
Alexander polynomial of $K$ associated to $\widehat{\rho}$. 
Then we prove that for any 2-bridge knot $K(r)$ in $H(p)$, 
$\widetilde{\Delta}_{\widehat{\rho},K}(t)$ is of the form 
$\left\{\dfrac{{\Delta}_{K(r)} (t)}{1-t}\right\} 
f(t) f(-t)$ for
some integer polynomial $f(t)$, where $H(p)$ is the set of $2$-bridge knots
$K(r), 0<r<1$, such that $G(K(r))$ is mapped onto a non-trivial free
product $\ZZ/2 * \ZZ/p$. Further, it is proved that 
$f(t) \equiv \left\{\dfrac{{\Delta}_{K} (t)} {1+t}\right\}^n$ (mod $p$), 
where ${\Delta}_{K} (t)$ is the Alexander polynomial of $K$.
Later we discuss the twisted Alexander polynomial 
associated to the general metacyclic representation.
\end{abstract}



\section{Introduction}

In the previous paper \cite{HM}, 
we studied the parabolic representation of the
group of a $2$-bridge knot and showed some properties of its twisted
Alexander polynomial. 
In this paper, we consider a metacyclic
representations of the knot group.

Let $G(m,p|k)$ be a (non-abelian) semi-direct product of two cyclic groups
$\ZZ/m$ and $\ZZ/p$, $p$ an odd prime, 
with the following presentation:
\begin{equation}
G(m,p|k)=\langle s,a|s^m=a^p=1,sas^{-1}=a^k\rangle,
\end{equation}
where $k$ is a primitive $m$-th root of 1 (mod $p$), i.e. 
$k^m \equiv 1$ (mod $p$), but $k^q \not\equiv 1$ (mod $p$) for any 
$q, 0<q<m$ and $k \ne 0, 1$.

If $k=-1$, then $m=2$ and hence $G(2,p|-1)$ is a dihedral group $D_p$.
Since $k$ is a primitive $m$-th root of 1 (mod $p$), 
$G(m,p|k)$ is imbedded in the symmetric
group $S_p$ and hence in $GL(p,\ZZ)$ via permutation matrices.

Now suppose that the knot group $G(K)$ of a knot $K$ is mapped onto
$G(m,p|k)$ for some $m,p$ and $k$. 
Then, we have a representation $f:\ G(K) 
\rightarrow G(m,p|k) \rightarrow GL(p,\ZZ)$
and the twisted Alexander polynomial 
$\widetilde{\Delta}_{f,K}(t)$
associated to $f$ is defined 
\cite{L} \cite{W} \cite{KL}.
One of our objectives is to
characterize these twisted Alexander polynomials.
In fact, we propose the following conjecture.

\medskip
\noindent{\bf Conjecture A}.
{\it
$\widetilde{\Delta}_{f,K}(t)=
\left\{\dfrac{\Delta_{K} (t)}{1-t}\right\} F(t)$,
where $\Delta_{K}(t)$ is the Alexander polynomial of $K$ and $F(t)$ is an
integer polynomial in $t^m$.
}
\medskip

First we study the case $k=-1$, dihedral representations of the knot group.
Let $D_p$ be a dihedral group of order $2p$, 
where $p=2n+1$ and $p$ is a prime. 
Then the knot group $G(K)$ of a knot $K$ is mapped onto $D_p$
if and only if ${\Delta}_{K} (-1) \equiv 0$ (mod $p$) 
\cite{Fox62}, \cite{Ha}.
Therefore, if ${\Delta}_{K} (-1) \ne \pm 1$,
$G(K)$ has at least one representation on a certain dihedral group $D_p$.
For these cases, we can make Conjecture A slightly sharper:

\medskip
\noindent{\bf Conjecture B}.
{\it
Let $\widehat{\rho}:\ 
G(K) \rightarrow D_p \rightarrow 
GL(p, \ZZ)$ be a non-abelian
representation of the knot group $G(K)$ of a knot $K$ and let 
$\widetilde{\Delta}_{\widehat{\rho}, K} (t)$ be the
twisted Alexander polynomial of $K$ associated to $\widehat{\rho}$. 
Then
\begin{equation}
\widetilde{\Delta}_{\widehat{\rho}, K} (t) =
\left\{\frac{{\Delta}_{K}(t)}{1-t}\right\}f(t) f(-t),
\end{equation}
where $f(t)$ is an integer polynomial and further,
\begin{equation}
f(t) \equiv \left\{\frac{{\Delta}_{K} (t)}{1+t}\right\}^{n} 
({\rm mod\ } p)
\end{equation}
}
\medskip

We should note that $(1+t)^2$ divides $\Delta_K(t)$ (mod $p$)
if and only if $\Delta_K(-1)\equiv 0$ (mod $p$).

The main purpose of this paper is to prove (1.2) for a 
$2$-bridge knot $K(r)$
in $H(p)$, $p$ a prime, and (1.3) for a $2$-bridge knot with
${\Delta}_{K}(-1) \equiv 0$ (mod $p$). 
(See Theorem 2.2.)
Here $H(p)$ is the set of $2$-bridge knots $K(r), 0<r<1,$ 
such that $G(K(r))$ is mapped onto a free product $\ZZ/2 * \ZZ/p$.
We note that knots in $H(p)$ have been studied extensively in \cite{GR} and \cite{ORS}.

A proof of the main theorem (Theorem 2.2) is given in 
Section 2 through
Section 7.
Since this paper is a sequel of \cite{HM}, 
we occasionally skip some details if
the argument used in \cite{HM} also works in this paper.

In Section 8, we consider another type of metacyclic groups, denoted by
$N(q,p)$.
$N(q,p)$ is a semi-direct product of two cyclic groups,
$\ZZ/2q$ and 
$\ZZ/p$
defined by
\begin{equation}
N(q,p)=\langle s, a| s^{2q}=a^p=1, sas^{-1}=a^{-1}\rangle,
\end{equation}
where $q \geq 1$ and $p$ is an odd prime and $\gcd(q,p)=1$. We note that
$N(1,p) =D_p$ and $N(2,p)$ is called a binary dihedral group.

Let $\widetilde{\nu}:\ G(K) \longrightarrow N(q,p) \longrightarrow
GL(2pq,\ZZ)$ be a representation of $G(K)$. 
(For details, see Section 8.)
Then we show that for a $2$-bridge knot $K(r)$, the twisted Alexander
polynomial $\widetilde{\Delta}_{\widetilde{\nu},K(r)}(t)$
associated to $\widetilde{\nu}$ is completely determined by
the Alexander polynomial
$\Delta_{K(r)}(t)$ and the twisted Alexander polynomial
$\widetilde{\Delta}_{\widehat{\rho},K(r)}(t)$ associated to 
$\widehat{\rho}$. (Proposition \ref{prop:8.5})

In Section 9, we give examples that illustrate our main 
theorem and Proposition \ref{prop:8.5}. 
It is interesting to observe that
$\widetilde{\Delta}_{\widetilde{\nu},K(r)}(t)$ 
is an integer polynomial in $t^{2q}$.
In Section 10, we briefly discuss general 
$G(m,p|k)$-representations of 
the knot group and give several examples, 
one of which is not a $2$-bridge knot,
that support Conjecture A.
In Section 11, we prove Proposition 2.1 and Lemma 5.2 
that plays a key role
in our proof of the main theorem.

Finally, for convenience, we draw a diagram below consisting of homomorphisms that
connect various groups and rings.

\begin{center}
$
\begin{array}{
cccc
ccc}
 & & GL(p,\ZZ) & & GL(2n,\ZZ)& &  \\
 & & \mbox{\large $\pi$}\uparrow& \nearrow& \hspace*{-11mm}
 \mbox{\large $\pi_0$}
  \ \ \ \ \ \uparrow \mbox{\large $\gamma$}& &  \\
G(K)&\underset{\mbox{\large $\rho$}}{\longrightarrow}&D_p&
\underset{\mbox{\large $\xi$}}{\longrightarrow}&GL(2,\CC)& &  \\
\downarrow& & \downarrow& & & &  \\
\ZZ G(K)&\longrightarrow&\ZZ D_p&
\underset{\mbox{\large $\zeta$}}{\longrightarrow}&\widetilde{A}(\omega)& & 
M_{2n,2n}(\ZZ[t^{\pm 1}])\\
& \mbox{\large $\rho^{*}$}\hspace{-1mm}\searrow&\downarrow& &\downarrow&&\uparrow 
\mbox{\large $\gamma^{*}$}\\
& & \ZZ D_p[t^{\pm1}]&\underset{\mbox{\large $\zeta^{*}$}}{\longrightarrow}&
\widetilde{A}(\omega)[t^{\pm 1}]&
\underset{\mbox{\large $\xi^{*}$}}{\longrightarrow}& 
M_{2,2}\bigl((\ZZ [\omega])[t^{\pm 1}]\bigr)
 \end{array}
 $
 \end{center}
 
 Here, $\tau=\rho\circ\xi, \widehat\rho=\rho\circ\pi, \rho_0=\rho\circ\pi_0,
 \eta=\xi\circ\gamma,
 \Phi^*=\rho^*\circ\zeta^*\circ\xi^*$ and $\nu=\rho\circ\xi\circ\gamma$.
 Unmarked arrows indicate natural extensions of homomorphisms.
 

\section{
Dihedral representations and statement of the main theorem}

We begin with a precise formulation of representations.
Let $p=2n+1$ and $D_p$ be a dihedral group of order $2p$ 
with a presentation:
$D_p = \langle x,y| x^2 = y^2 = (xy)^p = 1\rangle$.
As is well known, $D_p$ can be faithfully represented 
in $GL(p,\ZZ)$ by the map $\pi$ defined by:
\begin{center}$
x\mapsto
\left[
\begin{array}{ccccccc}
1&0&0&\cdots&0&0&0\\
0&0&0&\cdots&0&0&1\\
0&0&0&\cdots&0&1&0\\
\vdots & \vdots& \vdots& &\nana & &\vdots \\
\vdots &\vdots & &\nana& & \vdots& \vdots\\
\vdots  & & \nana & & \vdots&\vdots & \vdots\\
0&1&0&\cdots&0&0&0
\end{array}
\right]
\ 
y\mapsto
\left[
\begin{array}{ccccccc}
0&1&0&\cdots&0&0&0\\
1&0&0&\cdots&0&0&0\\
0&0&0&\cdots&0&0&1\\
0&0&0&\cdots&0&1&0\\
\vdots & \vdots& \vdots& & \nana& & \vdots\\
\vdots  & \vdots&  &\nana & \vdots& \vdots&\vdots \\
0&0&1&\cdots&0&0&0
\end{array}
\right]
$
\end{center}

However, $\pi$ is reducible. 
In fact, $\pi$ is equivalent to $id \ast{\pi}_0$, where
\begin{equation}
\pi_0: x\mapsto
\left[
\begin{array}{cccccc}
0&0&\cdots&0&0&1\\
0&0&\cdots&0&1&0\\
\vdots & \vdots& &\nana & & \vdots\\
\vdots& &\nana  & & \vdots& \vdots\\
0&1&\cdots&0&0&0\\
1&0&\cdots&0&0&0
\end{array}
\right]
y\mapsto
\left[
\begin{array}{ccccccc}
-1&0&0&\cdots&0&0&0\\
-1&0&0&\cdots&0&0&1\\
-1&0&0&\cdots&0&1&0\\
\vdots  & \vdots& \vdots& &  \nana&& \vdots\\
\vdots  & \vdots& & \nana & & \vdots& \vdots\\
  -1&0&1&\cdots&0&0&0\\
-1&1&0&\cdots&0&0&0
\end{array}
\right]
\end{equation}

For convenience, ${\pi}_0$ is called the
{\it irreducible representation} of $D_p$ (of degree $p-1=2n$).

Now let $K(r), 0<r<1, r= \frac{\beta}{\alpha}$ and $\gcd(\alpha,\beta)=1$, 
be a $2$-bridge knot and consider a Wirtinger presentation of
the group $G(K(r))$:
\begin{align}
&G(K(r)) = \langle x,y| R \rangle,\
{\rm where}\nonumber\\ 
&R = WxW ^{-1} y^{-1},
W = x^{\epsilon_1} 
y^{\epsilon_2} \cdots x^{\epsilon_{\alpha-2}} y^{\epsilon_{\alpha-1}}\ 
{\rm and}\nonumber\\
&\epsilon_j = \pm 1\ {\rm  for}\ 
1\leq j \leq \alpha -1.
\end{align}

Suppose $p$ be a prime. If $\alpha \equiv 0$ (mod $p$), then a mapping
\begin{equation}
\rho: x \mapsto x\ {\rm  and}\ y \mapsto y
\end{equation}
defines a surjection from $G(K(r))$ to $D_p$.

Therefore $\rho_0 = \rho \circ \pi_0$ defines a representation of
$G(K(r))$ into $GL(2n,\ZZ)$ and we can define the twisted Alexander
polynomial $\widetilde{\Delta}_{\rho_0, K(r)} (t)$ associated
to $\rho_0$. 
Since $\pi = id \ast \pi_0$, the twisted Alexander polynomial
associated to $\widehat{\rho} = \rho \circ \pi $ is given by
$\left[\dfrac{\Delta_{K(r)}(t)}{1-t}\right] \widetilde{\Delta}_{\rho_0, K(r)}(t)$ and
hence (1.2) becomes
\begin{equation}
\widetilde{\Delta}_{\rho_0, K(r)} (t) = f(t) f(-t).
\end{equation}

Now there is another representation of $D_p$ in $GL(2,\CC)$. To be more
precise, consider $\xi: D_p \rightarrow GL(2,\CC)$ given by
\begin{equation}
\xi (x) =\mtx{-1}{1}{0}{1}\ {\rm  and}\ \xi (y) 
=\mtx{-1}{0}{\omega}{1},
\end{equation}
where $\omega \in \CC$ is determined as follows.

First we set $\xi (x) = \mtx{-1}{1}{0}{1}$
and $\xi (y) =\mtx{-1}{0}{z}{1}$, and write
$\xi((xy)^k) = \mtx{a_k(z)}{b_k(z)}{c_k(z)}{d_k(z)}$.
Since $\xi(xy) = \mtx{1+z}{1}{z}{1}$,
we see that $a_k ,b_k ,c_k$
and $d_k$ are exactly the same polynomials found in \cite[(4.1)]{HM}. 
Further, as is
mentioned in \cite{HM}, $a_n(z)$ and $b_n(z)$ are given as follows: 
\cite[Propositions 10.2 and 2.4]{HM}:
\begin{equation}
a_n(z) =  \sum_{k=0}^{n}  \binom{n+k}{2k}z^k\ 
{\rm and}\
b_n(z)=\sum_{k=0}^{n-1} 
\binom{n+k}{2k+1} z^k.
\end{equation}

Since $(xy)^{2n+1} =1$, we have $(xy)^n x= y (xy)^n$ and
hence, a simple calculation shows that 
$\xi ((xy)^n x)=\xi (y (xy)^n)$ yields $a_n(z) + 2 b_n(z)= 0$. 
Therefore, the number $\omega$ we are looking for is a root
of $\theta_n (z) = a_n (z) + 2 b_n (z)$. Write
$\theta_n (z)={c_0}^{(n)} + {c_1}^{(n)} z + \cdots +
{c_{n-1}}^{(n)} z^{n-1} + {c_n}^{(n)} z^n$. 
Then we see
\begin{equation}
{\displaystyle
{c_k}^{(n)}=\binom{n+k}{2k}+ 2\binom{n+k}{2k+1}
=\frac{2n+1}{2k+1}\binom{n+k}{n-k}}.
\end{equation}

If $p = 2n+1$ is prime, then, for $0 \leq k \leq n-1$, 
${c_k}^{(n)} \equiv 0$
(mod $p$), but ${c_0}^{(n)} = p$ and ${c_n}^{(n)} = 1$. 
Therefore, by Eisenstein's criterion, 
$\theta_n (z)$ is irreducible and it is the minimal
polynomial of $\omega$.

Let $C_n$ be the companion matrix of $\theta_n (z)$. By substituting
$C_n$ for $\omega$, we have a homomorphism $\gamma : GL(2,\CC)
\rightarrow GL(2n,\ZZ)$, namely, $\gamma (1) = E_n$ and 
$\gamma (\omega)= C_n$, where $E_n$ is the identity matrix, and hence we obtain another
representation $\eta=\xi \circ \gamma : D_p \rightarrow GL(2n,\ZZ)$.

The following proposition is likely known, but since we are unable to find a
reference, we prove it in Section 11.

\begin{prop}\label{prop:2.1}
Two representations $\pi_0$ and $\eta$ are equivalent.
In other words, there is a matrix $U_n \in GL(2n,\ZZ)$ such that
\begin{equation}
U_n \pi_0 (x) {U_n}^{-1}= \eta (x)\ {\rm and}\ 
U_n \pi_0 (y) 
{U_n}^{-1} = \eta (y).
\end{equation}
\end{prop}

Let $K(r)$ be a $2$-bridge knot in $H(p)$. 
Then $\tau = \rho \circ \xi : 
G(K(r)) \rightarrow D_p \rightarrow GL(2,\CC)$ defines a
representation of $G(K(r))$ and let 
$\widetilde{\Delta}_{\tau, K(r)} (t |\omega)$ be
the twisted Alexander polynomial associated to $\tau$. 
Sometimes, we use the notation
$\widetilde{\Delta}_{\tau, K(r)} (t|\omega)$ to emphasize that the
polynomial involves $\omega$. 
Let $\omega_1, \omega_2, \cdots, \omega_n$ 
be all the roots of $\theta_n (t)$. 
Since $\theta_n (t)$ is
irreducible, the total $\tau$-twisted Alexander polynomial 
$D_{\tau, K(r)} (t)$ defined in \cite{SW} is given by
\begin{equation}
D_{\tau, K(r)} (t)=\prod_{j=1}^n
\widetilde{\Delta}_{\tau, K(r)} (t|\omega_j).
\end{equation}

It is known that the polynomial $D_{\tau, K(r)} (t)$ is rewritten as
\begin{equation}
D_{\tau, K(r)} (t)=\det[ \widetilde{\Delta}_{\tau,
K(r)} (t|\omega)]^{\gamma}.
\end{equation}

By (2.5), we see that $D_{\tau, K(r)} (t)$ is exactly the twisted
Alexander polynomial of $K(r)$ associated to $\nu = \rho \circ \eta
:G(K) \rightarrow GL(2n,\ZZ)$. 
Since, by Proposition \ref{prop:2.1}, 
$\pi_0$ and $\eta$ are equivalent, 
$\rho_0$ and $\nu$ are equivalent, and
hence $\widetilde{\Delta}_{\rho_0, K(r)} (t) = 
D_{\tau, K(r)}(t)$.

Conjecture A now becomes the following theorem under our
assumptions that will be proven in Sections 5-7.

\begin{thm}\label{thm:2.2}
If a 2-bridge knot $K(r)$ is in $H(p)$, then
\begin{equation}
D_{\tau, K(r)} (t) = f(t) f(-t)
\end{equation}
for some integer polynomial $f(t)$, 
and further, for any $2$-bridge knot
$K(r)$ with $\Delta_{K(r)}(-1) \equiv 0$ (mod $p$),
\begin{align}
&(1)\ D_{\tau, K(r)} (t) \equiv f(t) f(-t)\  {\rm (mod}\ p)\ 
{\rm and}\nonumber\\
&(2)\ f(t) \equiv \Bigl\{\dfrac{\Delta_K(t)}{1+t}\Bigr\}^n\ 
{\rm (mod}\ p),
\end{align}
where $\Delta_{K(r)}(t)$ is the Alexander polynomial of $K(r)$.
\end{thm}

We note that $\Delta_{K(r)}(t)$ is divisible by $1+t$ in $(\ZZ/p)[t^{\pm 1}]$.

\begin{rem}
If $n=1$, i.e., $p=3$, $\theta_1(z)=z+3$, and hence $\omega=-3$.
Therefore, $\gamma$ is an identity homomorphism and
$\widetilde\Delta_{\rho_0,K(r)}(t)=D_{\tau,K(r)}(t)$.
\end{rem}

\section{Basic formulas}
In this section, we list various formulas involving 
$a_k, b_k, c_k$ and $d_k $
which will be used throughout this paper. 
Most of these materials are
collected from Section 4 in \cite{HM}.

For simplicity, let 
$\xi (x) =X= \mtx{-1}{1}{0}{1}$ and
$\xi (y)=Y = \mtx{-1}{0}{\omega}{1}$,
where $\omega$ is a root of $\theta_n (z)$.

First we list several formulas which are similar to \cite[Proposition 4.2]{HM}

\begin{prop}\label{prop:3.1}
As before, write $(XY)^k$ = $
\mtx{a_k}{b_k}{c_k}{d_k}$.

\begin{align}
&(I)\  a_0 = d_0 = 1\ {\it and}\ b_0 = c_0 = 0.\nonumber\\
&(II)\  a_1 = 1+ \omega, b_1 =1, c_1 = \omega\  {\it and}\  d_1 = 1.\nonumber\\
&(III)\  (i)\ {\it For}\ k \geq 2,\nonumber\\
&\ \ (1)\ a_k = (2+ \omega) a_{k-1} - a_{k-2},\nonumber\\
&\ \ (2)\ \omega b_k = (1+ \omega) a_{k-1} - a_{k-2},\nonumber\\
&\ \ \ \ \ \ \ (ii)\  {\it For}\ k \geq 1,\nonumber\\
&\ \ (3)\ \omega b_k = a_k - a_{k-1},\nonumber\\
&\ \ (4)\ \omega b_k = c_k,\nonumber\\
&\ \ (5)\ a_k = \omega b_k + d_k,\nonumber\\
&\ \ (6)\ d_k = a_{k-1},\nonumber\\
&\ \ (7)\ b_k = b_{k-1} + a_{k-1},\nonumber\\
&\ \ (8)\ c_k + d_k = a_k,\nonumber\\
&\ \ (9)\ a_0 + a_1 + \cdots + a_{k-1} = b_k.
\end{align}
\end{prop}

Since a proof of Proposition \ref{prop:3.1} 
is exactly the same as that of Proposition 4.2 in \cite{HM}, we omit the details.

Next three propositions are different from the corresponding proposition
\cite[Proposition 4.4]{HM}, since they depend on defining relations of $D_p$.

\begin{prop}\label{prop:3.2}
Let $p=2n+1$.
\begin{align}
&(1)\ {\it For}\ 
0 \leq k \leq 2n, a_k = a_{2n-k}\ {\it and}\ 
a_{2n+1} = a_0.\nonumber\\
&(2)\ {\it For}\ 0 \leq k \leq 2n, 
b_k = -b_{p-k}\ {\it and}\ b_p =0.
\end{align}
\end{prop}

{\it Proof.}
Since $(XY)^k = (YX)^{p-k} = Y(XY)^{p-k}Y$, we have\\
$\mtx{a_k}{b_k}{c_k}{d_k}=\mtx{a_{p-k}-\omega b_{p-k}}{-b_{p-k}}
{-\omega a_{p-k}-c_{p-k}+\omega^2 b_{p-k}+\omega d_{p-k}}
{\omega b_{p-k}+d_{p-k}}$
and hence $a_k = a_{p-k} - \omega b_{p-k}$ and $b_k = - b_{p-k}$ which
proves (2). Further, $a_k = a_{p-k} - \omega b_{p-k} = a_{p-k} + \omega
b_k$ and thus, $a_{p-k} = a_{k-1}$ by (3.1)(III)(3). This proves (1). 
Finally, it is obvious that $a_p = a_0$.
\qed

\begin{prop}\label{prop:3.3}
Let $p=2n+1$. Then we have the following
\begin{align} 
&(1)\ a_0 + a_1 + \cdots + a_{2n} = 0,\nonumber\\
&(2)\ b_1 + b_2 + \cdots + b_{2n} = 0,\nonumber\\
&(3)\ d_0 + d_1 + \cdots + d_{2n} = 0,\nonumber\\
&(4)\ a_n + 2b_n =0.\nonumber\\
&(5)\ {\it If}\ 
k \equiv \ell\ {\rm (mod}\ p), \ 
{\it then}\ a_k = a_{\ell}, b_k =
b_{\ell}, c_k = c_{\ell}\ {\it and}\ d_k = d_{\ell}.
\end{align}
\end{prop}

{\it Proof.} First, we see that $(XY)^n X= Y(XY)^n$ implies\\
$\mtx{-a_n}{a_n+b_n}{-c_n}{c_n+d_n}=
\mtx{-a_n}{-b_n}{\omega a_n+c_n}{\omega b_n+d_n}$,
and hence
$a_n + b_n = - b_n$ that proves (4). (5) is immediate, since
$(XY)^p$ = 1. (1) follows from (3.1)(III)(9), since $a_0 + a_1 +
\cdots + a_{2n}=b_{2n+1}=0$. To show (2), use (3.1)(III)(3). Since
$b_0$ = 0, we see $\omega (b_1 + b_2 + \cdots + b_{2n})
=(a_1 - a_0) + (a_2 - a_1) + \cdots + (a_{2n-1} - a_{2n-2}) + (a_{2n} -
a_{2n-1}) 
=a_{2n} - a_0$ = 0, by (3.2)(1). (3) follows from
(3.1)(III)(6), since $d_0 = 1 = a_0 = a_{2n}$ and $d_0 + d_1 + \cdots
+ d_{2n} 
=1 + a_0 + a_1 + \cdots + a_{2n-1}$ = $a_0 + a_1 + \cdots +
a_{2n-1} + a_{2n} = 0$. 
\qed

Now we define an algebra $\widetilde {A}(\omega)$ using the group ring
$\ZZ D_p$.
Consider the linear extension $\widehat {\xi}$ of $\xi: \ZZ D_p
\rightarrow M_{2,2}(\ZZ[\omega])$ given by $\widehat {\xi} (x)=X$ and
$\widehat {\xi} (y)=Y$, where $M_{k,k} (R)$ denotes the ring of 
$k \times k$ matrices over a commutative ring $R$. 
Let ${\widehat {\xi}}^{-1} (0)$
be the kernel of $\widehat {\xi}$. 
Then $\widetilde {A}(\omega)=
\ZZ D_p / {\widehat {\xi}}^{-1} (0)$ is a non-commutative
$\ZZ[\omega]$-algebra. 
Some elements of ${\widehat {\xi}}^{-1} (0)$ can be
found in Proposition \ref{prop:3.4} 
below.

We define $\zeta: \ZZ D_p \rightarrow \widetilde {A}(\omega)$ 
to be the natural projection.

\begin{prop}\label{prop:3.4}
In $\widetilde{A}(\omega)$, the following formulas hold,
where $1$ denotes the identity of $\widetilde{A}(\omega)$.
\begin{equation}
{\it For}\ 1 \leq k \leq n, 
(xy)^k + (yx)^k = ( a_{k-1} + a_k) 1.
\end{equation}
\begin{align}
&(1)\ {\it For}\ 
1 \leq k \leq n-1, (xy)^k x + y(xy)^k = a_k (x+y),\nonumber\\
&(2)\ (xy)^n x = y(xy)^n=\frac{a_n}{2} (x+y) = -b_n (x+y).
\end{align}
\end{prop}

{\it Proof.}
To prove (3.4), it suffices to show that 
$(XY)^k + (YX)^k =(a_{k-1}+ a_k) E_n$. 
In fact, for $1 \leq k \leq n$, 
\begin{equation*}
(XY)^k + (YX)^k =(XY)^k + (XY)^{p-k}
=
\mtx{a_k+a_{p-k}}{b_k+b_{p-k}}{c_k+c_{p-k}}{d_k+d_{p-k}}.
\end{equation*}
Since $a_k + a_{p-k} = a_k + a_{k-1}$ by (3.2)(1), $b_k + b_{p-k}=0$ by
(3.2)(2), $c_k + c_{p-k} = \omega (b_k + b_{p-k}) =0$ and 
$d_k +d_{p-k} = a_{k-1} + a_{2n-k} = a_{k-1} + a_k$ by (3.1)(6) and (3.2)(1),
(3.4) follows immediately.
Next, for $1 \leq k \leq n-1$,
$(XY)^kX + Y(XY)^k= \mtx{-2a_k}{a_k}{\omega a_k}{2(c_k+d_k)}=a_k(X+Y)$, 
which proves (3.5)(1). Finally, (3.5)(2) follows, since $(xy)^n x
=y(xy)^n$ and $a_n = -2b_n$. 
\qed

\section{Polynomials over $\widetilde{A}(\omega)$}

In this section, as the first step toward a proof of 
Theorem \ref{thm:2.2}, 
we introduce one of our key concepts in this paper.

\begin{dfn}
 Let $\varphi (t)$ be a polynomial on $t^{\pm 1}$ with
coefficients in the non-commutative algebra $\widetilde{A}(\omega)$. We
say $\varphi (t)$ is {\it split} if $\varphi (t)$ is of the form:\\
$\varphi (t)=\sum_j \alpha_j t^{2j} + \sum_k \beta_k (x+y)
t^{2k+1}$, where $\alpha_j, \beta_k \in \ZZ[\omega]$.
The set of split polynomials is denoted by $S(t)$. 
For example, $\varphi(t)$ = $1+t^2, (x+y)t$ are split.
\end{dfn}

First we show that $S(t)$ is a commutative ring.

\begin{prop}\label{prop:4.2}
If $\varphi (t)$ and ${\varphi}^{\prime}(t)$ are split,
so are $\varphi (t) + {\varphi}^{\prime}(t)$ and $\varphi (t)
{\varphi}^{\prime}(t)$.
\end{prop}

{\it Proof.} Let 
$\varphi (t)=\sum_j \alpha_j t^{2j} + \sum_k \beta_k (x+y)
t^{2k+1}$ and ${\varphi}^{\prime}(t)=\sum_{\ell}
{\alpha_{\ell}}^{\prime}t^{2\ell} + \sum_m {\beta_m}^{\prime}(x+y)
t^{2m+1}$. 
Then obviously $\varphi (t)+ {\varphi}^{\prime}(t)$ is split. 
Further, 
\begin{align*}
\varphi (t) {\varphi}^{\prime}(t)& = \sum_{j, \ell} \alpha_j
{\alpha_{\ell}}^{\prime}t^{2j+2\ell} + \sum_{j,m} \alpha_j
{\beta_m}^{\prime}(x+y) t^{2j+2m+1}\\
&\ \ + \sum_{k,\ell} \beta_k {\alpha_\ell}^{\prime} (x+y) t^{2k+2\ell+1} +
\sum_{k,m} \beta_k {\beta_m}^{\prime}(x+y) (x+y) t^{2k+2m+2}.
\end{align*}
Since $(x+y)(x+y)=2+xy+yx=(2+b_2)1$ by (3.4) and (3.1)(III)(9), it follows
that $\varphi (t) {\varphi}^{\prime}(t)$ is split. 
\qed

Next, to obtain the proposition corresponding to Lemma 4.5 in \cite{HM}, we
define the polynomials over $\widetilde{A}(\omega)$. 

Let
$Q_k (t)=1 + (yx)t^2 + (yx)^2 t^4 + \cdots + (yx)^k t^{2k}$ and\\
$P_k (t)=1 + (xy)t^2 + (xy)^2 t^4 + \cdots + (xy)^k t^{2k}$. Note 
$Q_k (t) = y P_k (t) y$.
The following proposition is a slight modification of Lemma 4.5 in \cite{HM}.

\begin{prop}\label{prop:4.3}
Let $p=2n+1$.\\
(1) $(y^{-1} t^{-1}) (1-yt) Q_{2n} (t) yt(1-xt) \in S(t)$.\\
(2) $(y^{-1} t^{-1}) \{(1-yt) Q_n (t) yt + 
(yx)^{n+1}t^{2n+2}\} (1-xt)
\in S(t)$.\\
(3) $(y^{-1} t^{-1}) \{(1-yt) Q_{3n+1}(t) yt + 
(yx)^{3n+2}t^{6n+4}\}
(1-xt) \in S(t)$.\\
(4) $(y^{-1}t^{-1}) (1-yt) Q_{4n}(t) yt (1-xt) \in S(t)$.
\end{prop}

{\it Proof.}
 First we prove (2). Since
\begin{align*}
(1-yt) Q_n(t) yt + (yx)^{n+1}t^{2n+2}
&=(1-yt) yP_n(t) t + (yx)^{n+1} t^{2n+2}\\
&=yt (1-yt) P_n (t) + yt (xy)^n x t^{2n+1}\\
&=yt\{(1-yt)P_n (t) + (xy)^n x t^{2n+1}\}, 
\end{align*}
{\rm it\ suffices\ to\ show}
\begin{equation}
\{(1-yt) P_n (t) + (xy)^n x t^{2n+1}\}(1-xt) \in S(t).
\end{equation}

Now a simple computation shows that
\begin{align*}
&\{(1-yt) P_n (t) + (xy)^n x t^{2n+1}\}(1-xt)\\
&=\left\{\susum{k=0}{n} (xy)^k t^{2k}
- \susum{k=0}{n-1} y(xy)^k t^{2k+1}\right\} (1-xt)\\
&=1 + \susum{k=1}{n} \left\{(xy)^k
+ (yx)^k\right\} t^{2k} - \susum{k=0}{n-1} 
\left\{y(xy)^k + (xy)^k x\right\} t^{2k+1}\\
&=1+
\susum{k=1}{n} (a_{k-1} + a_k ) t^{2k} - 
\susum{k=0}{n-1} (x+y) a_k t^{2k+1}
\in S(t),
\end{align*} 
by (3.4) and (3.5).
This proves (4.1).

{\it Proof of (1).}
 Since 
 \begin{align*}(1-yt) Q_{2n}(t) yt (1-xt)
&=(1-yt) yP_{2n}(t) t
(1-xt)\\
&=yt (1-yt) P_{2n} (t) (1-xt),
\end{align*} it suffices to show
\begin{equation}
(1-yt) P_{2n} (t) (1-xt) \in S(t).
\end{equation}

However, the following straightforward 
calculation proves (4.2):
\begin{align*}
&(1-yt) P_{2n}(t) (1-xt)\\
&= {\textstyle\sum\limits_{k=0}^{2n}}(xy)^k t^{2k} -
{\textstyle\sum\limits_{k=0}^{2n}} y(xy)^k t^{2k+1} 
- {\textstyle\sum\limits_{k=0}^{2n}} (xy)^k x t^{2k+1} 
+ {\textstyle\sum\limits_{k=0}^{2n}} (yx)^{k+1} t^{2k+2}\\
&=1 + {\textstyle\sum\limits_{k=1}^{2n}} \left\{
(xy)^k + (yx)^k\right\} t^{2k} + (yx)^p t^{2p} 
- {\textstyle\sum\limits_{k=0}^{2n}}\left\{ y(xy)^k +
(xy)^k x\right\} t^{2k+1}\\
&=1+ {\textstyle\sum\limits_{k=1}^{2n}} (a_{k-1} + a_k) t^{2k} +
t^{2p} - {\textstyle\sum\limits_{k=0}^{2n}} 
a_k (x+y) t^{2k+1} \in S(t).
\end{align*}

{\it Proof of (3).}
 Since
\begin{align*}
&\left\{(1-yt) Q_{3n+1}(t) yt + (yx)^{3n+2} t^{6n+4}
\right\}(1-xt)\\
&=yt \{(1-yt)
P_{3n+1}(t) + (xy)^{3n+1} x t^{6n+3}\}(1-xt),
\end{align*} 
it suffices to show\\
\begin{equation}
\{(1-yt) P_{3n+1}(t) + (xy)^{3n+1}x t^{6n+3}\}
(1-xt) \in S(t).
\end{equation}
Since $P_{3n+1}(t)=P_{2n}(t) + t^{4n+2} P_n (t)$ and $(xy)^{3n+1} x =
(xy)^n x$, we must show\\
$\Bigl\{(1-yt) \{P_{2n}(t) + P_n (t) t^{4n+2}\} + 
(xy)^n x t^{6n+3}\Bigr\}(1-xt) \in 
S(t)$.
However, since $(1-yt) P_{2n}(t) (1-xt) \in S(t)$ by (4.2), it suffices
to show that
\begin{equation}
\{(1-yt) P_n (t) t^{4n+2} + (xy)^n x t^{6n+3}\}(1-xt) \in S(t).
\end{equation}
Now, (4.4) follows from (4.1), since $t^{4n+2}$ is split.

{\it Proof of (4).} 
Since $(yx)^{2n+1} = 1$, we have
\begin{equation*}
Q_{4n}(t)=\susum{k=0}{2n} 
(yx)^k t^{2k} + \susum{k=2n+1}{4n} (yx)^k t^{2k}
=(1+t^{2p})Q_{2n}(t). 
\end{equation*}
Since $(1+t^{2p})$ is split, it follows that
\begin{equation*}
(y^{-1}t^{-1}) (1-yt) Q_{4n}(t) yt (1-xt)
=(1+t^{2p})(y^{-1}t^{-1})
(1-yt) Q_{2n}(t) yt (1-xt)
\end{equation*} is split by (1). 
\qed

\section{Proof of Theorem 2.2.(I)}

In this section we prove Theorem \ref{thm:2.2} 
(2.11) for a torus knot $K(1/p)$, $p =2n+1$ a prime.
First we define various homomorphisms among group rings.\\
Let $g=x^{m_1} y^{m_2} x^{m_3}y^{m_4} \cdots x^{m_{k-1}} 
y^{m_k}$, where $m_j$ are integers and let $m = \susum{j=1}{k} 
m_j$ and $\ell$ is arbitrary. Then we have:
\begin{align}
&(1)\ {\rho}^\ast : \ZZ G(K) \rightarrow \ZZ D_p [t^{\pm 1}]\
{\rm is\ defined\ by}\ {\rho}^\ast (g)=\rho (g) t^m,\nonumber\\
&(2)\ {\zeta}^\ast:\ZZ D_p [t^{\pm 1}] \rightarrow
\widetilde{A}(\omega)[t^{\pm 1}]\
{\rm is\ defined\ by}\ {\zeta}^\ast
(gt^{\ell})=\zeta (g) t^{\ell},\nonumber\\
&(3)\ 
{\xi}^\ast: \widetilde{A}(\omega) [t^{\pm 1}] \rightarrow M_{2,2}
(\ZZ[\omega][t^{\pm 1}])\ {\rm is\ defined\ by}\ 
{\xi}^\ast (gt^{\ell})=\xi (g) t^{\ell},\nonumber\\
&(4)\ 
{\gamma}^\ast: M_{2,2}(\ZZ[\omega] [t^{\pm 1}]) \rightarrow M_{2n,2n} 
(\ZZ[t^{\pm 1}])\ {\rm is\ defined\ by} \nonumber\\
&\ \ \ {\gamma}^\ast 
\mtx{\sum_j p_j t^j}{\sum_j q_j t^j}
{\sum_j r_j t^j}{\sum_j s_j t^j}
=\mtx{\sum_j \gamma(p_j) t^j}{\sum_j \gamma(q_j) t^j}
{\sum_j \gamma(r_j) t^j}{\sum_j \gamma(s_j) t^j}.
\end{align}

Now we show the following proposition.

\begin{prop}\label{prop:5.1}
Let $p=2n+1$, a prime. Then $D_{\tau,K(1/p)}$(t) is of
the form $q(t) q(-t)$ for some integer polynomial $q(t)$.
\end{prop}

{\it Proof.} 
We write $G(K(1/p))= \langle x,y|
R_0 = W_0 x {W_0}^{-1} y^{-1} =1\rangle$, where
$W_0 = (xy)^n$.
Consider the free derivative of $R_0$ with respect to $x$;
\begin{equation*}
 \dfrac{\partial R_0}{\partial x}
=(1-y)
\dfrac{\partial W_0}{\partial x} + W_0
=(1-y) \susum{k=0}{n-1} 
(xy)^k + (xy)^n,
\end{equation*}
and we write
\begin{equation*}
{\Phi}^\ast\left( \dfrac{\partial R_0}{\partial x}\right)
=\mtx{h_{11}(t)}{h_{12}(t)}{h_{21}(t)}{h_{22}(t)},
\end{equation*}
where ${\Phi}^\ast=
{\rho}^\ast \circ {\zeta}^\ast \circ {\xi}^\ast$. \\
Then we see; 
\begin{align}
(1)\ h_{11}(t)&=\susum{k=0}{n} a_k t^{2k} + 
 \susum{k=0}{n-1} 
a_k t^{2k+1}
=(1+t) \susum{k=0}{n-1} a_k t^{2k} + a_n t^{2n},\nonumber\\
(2)\  h_{12}(t)&=\susum{k=0}{n} b_k t^{2k} + 
\susum{k=0}{n-1} 
b_k t^{2k+1}
=(1+t) \susum{k=0}{n-1} b_k t^{2k} + b_n t^{2n},\nonumber\\
(3)\ 
h_{21}(t)&=\susum{k=0}{n} c_k t^{2k} - \omega
\susum{k=0}{n-1} a_k t^{2k+1} - \susum{k=0}{n-1} 
c_k t^{2k+1}\nonumber\\
& =- \omega
t \susum{k=0}{n-1} a_k t^{2k} + (1-t) 
\susum{k=0}{n-1} c_k t^{2k} + c_n
t^{2n},\nonumber\\
(4)\  h_{22}(t)&=\susum{k=0}{n} d_k t^{2k} - \omega
\susum{k=0}{n-1} b_k t^{2k+1} - \susum{k=0}{n-1} d_k t^{2k+1}
\nonumber\\
&=- \omega t \susum{k=0}{n-1} b_k t^{2k} + (1-t) 
\susum{k=0}{n-1} d_k t^{2k} + d_nt^{2n}.
\end{align}

Since $h_{11}(1) = 0$ and $h_{21}(1) = 0$, both $h_{11} (t)$ and 
$h_{21}(t)$ are divisible by $1-t$.
 In fact, we have:
\begin{align*}
h_{11}(t)&=(1-t) \Bigl\{\ \ \susum{k=0}{n-1} 
(2a_0 + 2a_1 + \cdots +2a_{k-1} + a_k)t^{2k}\\
&\hspace*{18mm}
 + \susum{k=0}{n} (2a_0 + 2a_1 
+ \cdots + 2a_k)t^{2k+1}\Bigr\}\\
&
=(1-t) \left\{\susum{k=0}{n-1} (b_k + b_{k+1}) t^{2k} + 
\susum{k=0}{n-1} 
2b_{k+1} t^{2k+1}\right\},\ {\rm and}\\
h_{21}(t)&=- \omega t(1-t^2) 
\susum{k=0}{n-2}( a_0 + a_1 + \cdots
+a_k) t^{2k} \\
&\hspace*{15mm}- \omega t(1-t) 
(a_0 + a_1 + \cdots +a_{n-1}) t^{2n-2} 
+(1-t) \susum{k=1}{n-1} c_k t^{2k}\\
&
=(1-t)\left\{- \omega t(1+t)
\susum{k=0}{n-2} b_{k+1} t^{2k} - \omega tb_n t^{2n-2} + \susum{k=1}{n-1}
c_k t^{2k}\right\}.
\end{align*}

Since $c_k = \omega b_k$, we see
$h_{21}(t) = (1-t) \left\{ - \omega t 
\susum{k=0}{n-1}b_{k+1}t^{2k}\right\}$, and hence,\\
$\dfrac{1}{1-t} \det 
\left(\dfrac{\partial R_0}{\partial x}\right)^{{\Phi}^\ast}
=\dfrac{1}{1-t} \det \mtx{h_{11}(t)}{h_{12}(t)}
{h_{21}(t)}{h_{22}(t)}
=\det \mtx{h_{11}^{\prime}(t)}{ h_{12}(t)}
{h_{21}^{\prime}(t)}{h_{22}(t)}$,
where 
\begin{align*}
{h_{11}}^{\prime}(t)&=\susum{k=0}{n-1} 
(b_k + b_{k+1}) t^{2k} +
\susum{k=0}{n-1} 2b_{k+1} t^{2k+1}\\
&=\susum{k=1}{n-1} b_k (1+t^2)t^{2k-2} + 
b_nt^{2n-2}+\susum{k=1}{n} 2b_k t^{2k-1},\ {\rm and}\\
{h_{21}}^{\prime}(t)&=- \omega t \susum{k=0}{n-1} 
b_{k+1} t^{2k}.
\end{align*}

Let $g(t)= \susum{k=1}{n-1} b_k (1+t^2) t^{2k+2} 
+ b_n t^{2n-2}$ and
$h(t) = \susum{k=1}{n} b_k t^{2k-1}$.
Then 
\begin{equation*}
{h_{11}}^{\prime}(t) = g(t) +2h(t)\ {\rm  and}\
{h_{21}}^{\prime}(t) = - \omega h(t).
\end{equation*}
Further a straightforward computation shows that
\begin{equation*}
{h_{11}}^{\prime}(t) + h_{12}(t)=(1+t) (g(t) + h(t)).
\end{equation*}
And, 
\begin{align*}
{h_{21}}^{\prime}(t) + h_{22}(t)
&=- \omega t \susum{k=1}{n} b_k
t^{2k-2} - \omega t \susum{k=1}{n-1}b_k t^{2k} +(1-t) 
\susum{k=1}{n-1} 
d_k t^{2k} + d_n t^{2n}\\
&=- \susum{k=1}{n} c_k t^{2k-1} -
\susum{k=1}{n-1} c_k t^{2k+1} + (1-t) 
\susum{k=0}{n-1} d_k t^{2k} + d_n 
t^{2n}.
\end{align*}

Since $c_k + d_k = a_k$ and $d_0 = a_0$, we see
\begin{equation*}
- \susum{k=1}{n-1} c_k t^{2k+1} - 
\susum{k=0}{n-1}d_k t^{2k+1} =
- \susum{k=0}{n-1} a_k t^{2k+1},
\end{equation*}
and hence
\begin{equation*}
{h_{21}}^{\prime}(t) + h_{22}(t)=\susum{k=0}{n} 
d_k t^{2k} - 
\susum{k=0}{n-1}(a_k + c_{k+1}) t^{2k+1}.
\end{equation*}
Now, ${h_{21}}^{\prime}(t) + h_{22}(t)$ is divisible by 
$1+t$, and in
fact, we have
\begin{equation*}
{h_{21}}^{\prime}(t) + h_{22}(t)=(1+t)\{g(t) - 2h(t) 
- \omega h(t)\}.
\end{equation*}
Therefore, 
\begin{align*}
\dfrac{1}{(1-t)(1+t)} \det \left({\Phi}^\ast \dfrac{\partial
R_0}{\partial x}\right)
&=\det 
\mtx{g(t) + 2h(t)}{g(t) + h(t)}
{-\omega h(t)}{g(t) - 2h(t) - \omega h(t)}\\
&=\det 
\mtx{g(t) + 2h(t)}{- h(t)}
{- \omega h(t)}{ g(t) - 2h(t)},
\end{align*} 
and hence
\begin{equation}
\widetilde{\Delta}_{\tau,K(1/p)}(t) = g(t)^2 - (4+w) h(t)^2.
\end{equation}

Now we apply the following key lemma.

\begin{lemm}\label{lem:5.2}
Let $C_n$ be the companion matrix of $\theta_n (z)$, the
minimal polynomial of $\omega$. Then there exists a matrix $V_n \in
GL(n,\ZZ)$ such that
${V_n}^2 = 4 E_n + C_n$.
\end{lemm}

Since our proof involves a lot of computations, 
the proof is postponed to Section 11.

Since the total twisted Alexander polynomial of $K(1/p)$ at $\tau$ is
$D_{\tau,K(1/p)}(t) = \det
[\widetilde{\Delta}_{\tau,K(1/p)}(t)]^{{\gamma}^\ast}$, we obtain,
noting that $V_n$ commutes with $C_n$,
\begin{align*}
D_{\tau, K(1/p)}(t)&=\det[g(t|C_n)^2 - {V_n}^2 h(t|C_n)^2]\\
&=\det [g(t|C_n) - V_n h(t|C_n)] \det [g(t|C_n) +V_n h(t|C_n)].
\end{align*}
Let $q(t) = \det [g(t|C_n) - V_n h(t|C_n)]$. 
Then since $g(-t) = g(t)$
and $h(-t)= - h(t)$, it follows that 
\begin{equation*}
D_{\tau,K(1/p)}(t) = q(t) q(-t).
\end{equation*} 
This proves
Theorem \ref{thm:2.2}
(2.11) for $K(1/p)$.

\begin{rem}
It is quite likely that
\begin{equation}
q(t) = (1+t)^n \left\{\Delta_{K(1/p)}(t)
\right\}^{n-1},
\end{equation}
where $\Delta_{K(1/p)}(t)$ is the Alexander polynomial of $K(1/p)$.
\end{rem}

\section{Proof of Theorem 2.2 (II)}

Now we return to a proof of Theorem \ref{thm:2.2} 
(2.11) for a $2$-bridge knot $K(r)$ in $H(p)$.
Let $G(K(r)) = \langle x,y|R\rangle, R=WxW^{-1} y^{-1}$, 
be a Wirtinger presentation of $G(K(r))$.
Then as is shown in \cite{HM}, 
$R$ is written freely as a product of conjugates
of $R_0$:
$R=\prod_{j=1}^s u_j {R_0}^{\epsilon_j}{u_j}^{-1}$,
where for $1 \leq j \leq s$, $\epsilon_j = \pm 1$ and $u_j \in F(x,y)$,
the free group generated by $x$ and $y$, and 
$\frac{\partial R}{\partial x}=
\sum_{j} \epsilon_j u_j (\frac{\partial R_0}{\partial x})$, and
hence
\begin{align*}
 \widetilde{\Delta}_{\tau, K(r)}(t)
 &=\det\left(\frac{\partial R}{\partial
x}\right)^{{\Phi}^\ast}/ \det (y^{{\Phi}^{\ast}}- E_2)\\
&=
\widetilde{\Delta}_{\tau, K(1/p)}(t) \det \bigl(\sum_{j} \epsilon_j u_j\bigr)
^{{\Phi}^{\ast}}.
\end{align*}

As we did in \cite{HM}, we study $\lambda (r)=(\sum_{j} \epsilon_j u_j)
^{{\tau}^\ast}
\in \widetilde{A}(\omega)[t^{\pm 1}]$, 
where ${\tau}^\ast = 
{\rho}^{\ast} \circ {\zeta}^{\ast}$. 
For simplicity, we denote
${\tau}^{\ast}(\lambda (r))$ by 
${\lambda_r}^{\ast}(t)$. In fact, it is
a polynomial in $t^{\pm 1}$.\\
Since $K(r) \in H(p)$, the continued fraction of $r$ is of the form:\\
$r=[pk_1, 2m_1, pk_2, \cdots, 2m_{\ell}, pk_{\ell +1}]$, where $k_j$ and
$m_j$ are non-zero integers.\\
First we state the following proposition.

\begin{prop}\label{prop:6.1}   
Suppose $K(r)$ and $K(r^{\prime})$ belong to $H(p)$ and
let\\
 $r=[pk_1, 2m_1, pk_2, \cdots, 2m_{\ell}, pk_{\ell +1}]$, 
$r^{\prime}=[{pk_1}^{\prime}, 2m_1, {pk_2}^{\prime}, \cdots, 2m_{\ell},
pk^{\prime}_{\ell +1}]$ be continued fractions of $r$ and $r^{\prime}$. 
Suppose that $k_j \equiv {k_j}^{\prime}$ (mod $4$)
for each $j, 1 \leq j \leq \ell +1$. Then if $y^{-1}t^{-1}
\lambda_r^{\ast}(t)$ is split, so is $y^{-1} t^{-1} 
\lambda_{r^{\prime}}^{\ast}(t)$.
\end{prop}

Since a proof is analogous to that of Proposition 6.3 in \cite{HM}, 
we omit the
details.

Now we study the polynomial ${\lambda_r}^{\ast}(t) \in 
\widetilde{A}(\omega) [t^{\pm 1}]$ and we prove that 
$y^{-1}t^{-1}
{\lambda_r}^{\ast}(t)$ is split.
As is seen in Section 7 in \cite{HM}, 
${\lambda_r}^{\ast}(t)$ is written as
${w^\ast}_{2\ell +1}(t)$ and we will prove the following proposition. The
same
notation employed in Section 7 in \cite{HM} 
will be used in this section.

\begin{prop}\label{prop:6.2}
$ y^{-1} t^{-1}w^\ast_{2\ell +1}(t) \in S(t)$.
\end{prop}

{\it Proof.} Use induction on $j$. 
First we prove $y^{-1}t^{-1}{w^\ast}_1 (t) \in
S(t)$.\\
(1) If $w_1(t) = yt$, then $y^{-1}t^{-1}{w^\ast}_1 (t) = 1$ and hence
$y^{-1}t^{-1}{w^\ast}_1 (t) \in S(t)$.\\
(2) If $w_1 = y - (yx)^{n+1}$, then ${w^\ast}_1 (t) = yt - (yx)^{n+1}
t^{2n+2}$ and\\
 $y^{-1}t^{-1}{w^\ast}_1 (t)=1-(xy)^n x t^{2n+1}$ = $1
+ b_n (x+y) t^{2n+1}$ and hence\\
 $y^{-1}t^{-1}{w^\ast}_1 (t) \in S(t)$.\\
(3) If $w_1 = - (yx)^{n+1}$, then ${w^\ast}_1 (t) = - (yx)^{n+1}
t^{2n+2}$ and\\ $y^{-1}t^{-1}{w^\ast}_1 (t)
=-(xy)^n x t^{2n+1}=b_n
(x+y) t^{2n+1}$ and hence\\
 $y^{-1}t^{-1}{w^\ast}_1 (t) \in S(t)$.

Now suppose $y^{-1}t^{-1}{w^\ast}_{2j - 1} (t) \in S(t)$ for $j \leq
\ell$, and we claim \\
$y^{-1}t^{-1}{w^\ast}_{2\ell + 1}(t) \in S(t)$.
There are three cases to be considered. (See \cite[Proposition 7.1.]{HM}

Case 1. $k_{\ell+1}=1$.
$w_{2\ell + 1}
=
\{(1-y) Q_n y + (yx)^{n+1}\}
  \sum_{j} m_j (x-1)
y^{-1} w_{2j-1} - (yx)^{n+1} y^{-1} w_{2\ell-1} + y$.\\
Then 
\begin{align*}
y^{-1}t^{-1}{w^\ast}_{2\ell + 1}(t)
&=
y^{-1}t^{-1}\{(1-yt) Q_n 
(t) yt \\
&\ \ + (yx)^{n+1}t^{2n+2}\}
\sum_{j} m_j (xt-1)
y^{-1}t^{-1}{w^\ast}_{2j-1}(t)\\
&\ \ - (xy)^n x
t^{2n+1}(y^{-1}t^{-1}{w^\ast}_{2\ell - 1}(t)) + 1.
\end{align*} 
By Proposition \ref{prop:4.3}(2), 
each summand is split. 
Further, $-(xy)^n x t^{2n+1}= b_n (x+y) t^{2n+1} \in S(t)$ and 
$1 \in S(t)$. Therefore, the sum of them is split.

Proofs of the other cases are essentially the same.

Case 2. $k_{\ell+1}=2$.\\
$w_{2\ell+1}
=(1-y) Q_{2n} y\{\sum_{j}m_j (x-1)
y^{-1}w_{2j-1}\} + (yx)^{2n+1}w_{2\ell-1} - (yx)^{n+1} + y$.\\
Then 
\begin{align*}
y^{-1}t^{-1}{w^\ast}_{2\ell + 1}(t)
&=y^{-1} t^{-1}(1-yt)
Q_{2n}(t) yt \{ \sum_{j} m_j (xt-1)y^{-1}t^{-1}{w^\ast}_{2j - 1}(t)\}\\
&\ \ +y^{-1}t^{-1} t^{4n+2} {w^\ast}_{2\ell - 1}(t) - x(yx)^n t^{2n+1}+1.
\end{align*}

Again, $y^{-1}t^{-1}(1-yt) Q_{2n}(t) yt (xt-1) \in S(t)$ by Proposition
\ref{prop:4.3}(1)   
and\\ $y^{-1}t^{-1}{w^\ast}_{2j- 1}(t) \in S(t)$ by induction
hypothesis and $t^{4n+2}$,\\
 $- x(yx)^n t^{2n+1} = b_n (x+y) t^{2n+1}$ and
1 are split. Thus, $y^{-1}t^{-1}{w^\ast}_{2\ell + 1}(t) \in S(t)$.

Case 3. $k_{\ell+1}=3$.
\begin{align*}
w_{2\ell+1}
& = \{(1-y) Q_{3n+1}y + (yx)^{3n+2}\} 
\sum_{j} m_j
(x-1) y^{-1} w_{2j-1} \\
&\ \ - (yx)^{3n+2} y^{-1}w_{2\ell-1}+(yx)^{p}y - (yx)^{n+1}+y.
\end{align*}
Then 
\begin{align*}
y^{-1}t^{-1}{w^\ast}_{2\ell + 1}(t)
&=y^{-1} t^{-1}\{(1-yt)Q_{3n+1}(t) yt\\
&\ \  + (yx)^{3n+2}t^{6n+4}\} \sum_{j} m_j (xt-1)
y^{-1}t^{-1}{w^\ast}_{2j-1}(t)\\
&\ \  - (xy)^n x
t^{6n+3}(y^{-1}t^{-1}{w^\ast}_{2\ell - 1}(t))+t^{2p} - (xy)^n x t^{2n+1}+1.
\end{align*}

We see that $y^{-1}t^{-1}{w^\ast}_{2\ell + 1}(t)$ is split, since
each of 
$y^{-1}t^{-1}\{(1-yt) Q_{3n+1}(t) yt 
+ (yx)^{3n+2} t^{6n+4}\} (xt-1)$, 
$y^{-1}t^{-1}{w^\ast}_{2j - 1}(t)$ and
$- (xy)^n x t^{6n+3}=b_n (x+y) t^{6n+3}$
 and $- (xy)^n x t^{2n+1} = b_n (x+y)
t^{2n+1}$ is split. 
This proves Proposition \ref{prop:6.1}
\qed 

Now a proof of (2.11) for our knots is exactly the same as we did in Section 5.
Since $y^{-1}t^{-1}{w^\ast}_{2\ell + 1}(t) \in S(t)$, we can write
\begin{equation*}
y^{-1}t^{-1}{w^\ast}_{2\ell + 1}(t)
=\sum_{j} \alpha_j t^{2j} +
\sum_{k} \beta_k (x+y) t^{2k+1},
\end{equation*} 
where $\alpha_j, \beta_k \in \ZZ[\omega]$.

Define $g(t) = \sum_{j} \alpha_j t^{2j}$ and $h(t) = \sum_{k} \beta_k 
t^{2k+1}$. 
Since $X+Y = 
\mtx{-2}{1}{\omega}{2}$,
\begin{equation*}
{\xi}^\ast [y^{-1}t^{-1}{w^\ast}_{2\ell + 1}(t)]
=\mtx{g(t) - 2h(t)}{h(t)}
{\omega h(t)}{g(t)+2h(t)}\ 
{\rm and}
\end{equation*}
\begin{equation*}
\det (y^{-1}t^{-1}{w^\ast}_{2\ell + 1}
(t))^{{\xi}^\ast}
=
g(t)^2 - (\omega+4)
h(t)^2.
\end{equation*}

Thus $\widetilde{\Delta}_{\tau,K(r)}(t|\omega)
=\widetilde{\Delta}_{\tau, K(1/p)}(t|\omega) 
\bigl\{g(t)^2 - (w+4) h(t)^2\bigr\}$,
and hence, we have
\begin{equation*}
D_{\tau, K(r)}(t) 
= D_{\tau, K(1/p)}(t) \det[g(t|C_n)^2 - (C_n +4E_n)
h(t|C_n)^2].
\end{equation*}

Now by Lemma \ref{lem:5.2},
there exists a matrix $V_n \in GL(n, \ZZ)$ such
that ${V_n}^2 = C_n + 4E_n$. Since $V_n$ commutes with $C_n$, we see
\begin{align*}
g(t|C_n)^2 - (C_n +4E_n) h(t|C_n)^2
&=g(t|C_n)^2 - {V_n}^2 h(t|C_n)^2\\
&=\{g(t|C_n) - V_n h(t|C_n)\} \{g(t|C_n) + V_n h(t|C_n)\}.
\end{align*}

Let $f(t) = \det [g(t|C_n) - V_n h(t|C_n)]$. 
Since $h(-t|C_n) = -
h(t|C_n)$ and $g(- t|C_n)= g(t|C_n )$, $f(-t) = \det [g(t|C_n) + V_n 
h(t|C_n)]$, and thus,
\begin{equation*}
\det [g(t|C_n)^2 - (C_n +4E_n ) h(t|C_n)^2] = f(t) f(-t).
\end{equation*}

Therefore, $D_{\tau, K(r)}(t) = D_{\tau, K(1/p)}(t) f(t) f(-t)$. Since
$D_{\tau, K(1/p)}(t)$ is of the form $q(t) q(-t)$, it follows that
$D_{\tau, K(r)}(t) = F(t) F(-t)$, where $F(t)=q(t) f(t)$.

This proves (2.11) for $K(r)$ in $H(p)$.
\qed

\section{Proof of Theorem 2.2 (III)}

In this section, we prove (2.12) for a $2$-bridge knot $K(r)$ with
$\Delta_{K(r)} (-1) \equiv 0$ (mod $p$).

First we state the following easy lemma without proof.

\begin{lemm}\label{lem:7.1}
Let $M$ be a $2n \times 2n$ matrix over a commutative ring
which is
decomposed into four $n \times n$ matrices, $A, B, C$ and $D$:
$M =\mtx{A}{B}{C}{D}$.

Suppose that each matrix is lower triangular and in particular, $C$ is
strictly
lower triangular, namely, all diagonal entries are 0. Then
$\det M = (\det A) (\det D)$, and hence, $\det M$ is the product of all
diagonal entries of $M$.
\end{lemm}

Lemma \ref{lem:7.1} 
can be proven easily by induction on $n$.

Now let $K(r), 0<r<1$, be a $2$-bridge knot and 
consider a Wirtinger
presentation $G(K(r)) = \langle x,y| R\rangle$,
where $R=x^{\epsilon_1}y^{\eta_1} x^{\epsilon_2} y^{\eta_2} \cdots 
x^{\epsilon_{\alpha}} y^{\eta_{\alpha}}$ and $\epsilon_j, \eta_j = \pm 1$
for $1\leq j \leq \alpha$.

Applying the free differentiation, we have
$\dfrac{\partial R}{\partial x}
=\susum{i=1}{\alpha} g_i, g_i \in
\ZZ G(K)$,
where
\begin{equation}
g_i=
\begin{cases}\begin{array}{ll}
x^{\epsilon_1}y^{\eta_1} x^{\epsilon_2} y^{\eta_2}
\cdots x^{\epsilon_{i-1}} y^{\eta_{i-1}}& {\rm if\ } \epsilon_i =1\\
-x^{\epsilon_1}y^{\eta_1} x^{\epsilon_2} y^{\eta_2} \cdots 
x^{\epsilon_{i-1}} y^{\eta_{i-1}} x^{-1}& {\rm if\ } \epsilon_i = -1.
\end{array}
\end{cases}
\end{equation}
Let $\Psi: \ZZ G(K) \rightarrow \ZZ[t^{\pm 1}]$ be the homomorphism defined by
$ \Psi( g_i)=\epsilon_i t^{m_i}$, where $m_i = \susum{j=1}{i-1}
(\epsilon_j + \eta_j ) + \dfrac{\epsilon_i - 1}{2}$.

Then $\left(\dfrac{\partial R}{\partial x}\right)^\Psi$ gives the Alexander
polynomial $\Delta_{K(r)}(t)$ of $K(r)$. 
On the other hand,
$\dfrac{1}{(1-t)(1+t)} \det\left(
\dfrac{\partial R}{\partial x}\right)^{{\Phi}^\ast}$ gives
the twisted Alexander polynomial
$\widetilde{\Delta}_{\rho_0,K(r)}(t|w)$ associated to the irreducible
dihedral representation $\rho_0$, and
further, we see $D_{\tau,K(r)}(t)=\det\left[
\dfrac{1}{(1-t)(1+t)}
(\dfrac{\partial R}{\partial x})^{{\Phi}^\ast}\right]^{
{\gamma}^\ast}$.

Now using (7.1), we compute $\left(\dfrac{\partial R}{\partial
x}\right)^{
{\Phi}^\ast}=\sum_{i} {{\Phi}^\ast}(g_i)$.\\
If $\epsilon_i = 1$, then $m_i$ is even and
\begin{align*}
{{\Phi}^\ast}(g_i) &= 
[(xy)^{i-1}]^{\xi} t^{m_i}\\
&=\mtx{a_{i-1}}{b_{i-1}}{c_{i-1}}{d_{i-1}}t^{m_i}.
\end{align*}
If $\epsilon_i = - 1$, then
$m_i$ is odd and
\begin{align*}
{{\Phi}^\ast}(g_i) &= - [(xy)^{i-1} x]^{\xi} t^{m_i}\\
& =-\mtx{-a_{i-1}}{ a_{i-1} + b_{i-1}}
 {-c_{i-1}}{ c_{i-1} + d_{i-1}}t^{m_i}.
 \end{align*}

Therefore we have
\begin{align*}
(\dfrac{\partial R}{\partial x})^{{\Phi}^\ast}&=\susum{i}{}
{{\Phi}^\ast}(g_i)\\
&={\displaystyle \sum_{m_i = even}}
\mtx{a_{i-1}}{b_{i-1}}{c_{i-1}}{d_{i-1}} t^{m_i}
- {\displaystyle \sum_{m_j = odd}} 
\mtx{-a_{j-1}}{a_{j-1} + b_{j-1}}{-c_{j-1}}{c_{j-1} +
d_{j-1}}t^{m_j}.
\end{align*}
We note that as polynomials on $\omega$, the constant terms of $a_{i-1}$
and $d_{i-1}$ both are $1$. 
Further, since $c_{i-1}= \omega b_{i-1}$, the
constant term of $c_{i-1} + d_{i-1}$ is also 1, and hence
\begin{equation*}
\sum_{i} [{g_i}^{{\Phi}^\ast}]^{{\gamma}^\ast}
=\mtx{\Delta_{K(r)}(-t)
+ \omega \mu_{11}}{\mu_{12}}
{\omega \mu_{21}}
{\Delta_{K(r)}(t) +
\omega \mu_{11}},\
{\rm where}\ \mu_{ij} \in (\ZZ[\omega])[t^{\pm 1}].
\end{equation*}

If we replace $\ZZ$ by $\ZZ/p$, then $C_n$ is reduced to
$\left[\begin{array}{ccc|c}
0&\cdots&0&0\\ \hline
&&&0\\
&\Hsymb{E} & &\vdots\\
&&&0
\end{array}
\right]
$
and hence
$\sum_{i} [{g_i}^{\Phi^\ast}]^{\gamma^\ast} \equiv \mtx{A}{B}{C}{D}$ 
(mod $p$),
where $A,B,C$ and $D$ are lower triangular and in particular, $C$ is
strictly lower
triangular, and each diagonal entry of $A$ and $D$ is $\Delta_{K(r)}(t)$
(mod $p$) and $\Delta_{K(r)}(-t)$ (mod $p$), respectively. 
Therefore, by Lemma \ref{lem:7.1}, we have 
\begin{align*}
D_{\tau,K(r)}(t) &\equiv 
\det (\susum{i}{} [{g_i}^{\Phi^\ast}]^{{\gamma}^\ast}) / 
\det [(1-t)(1+t)]^{\gamma^\ast}\\
& \equiv 
\left\{\dfrac{\Delta_{K(r)}(t)}{1+t}\right\}^n 
\left\{
\dfrac{\Delta_{K(r)}(-t)}{1-t}\right\}^n\
{\rm (mod}\ p).
\end{align*}

This proves (2.12) for any $2$-bridge knot $K(r)$ with $\alpha \equiv 0$
(mod $p$).
We note that $\Delta_{K(r)}(t)$ is divisible by $1+t$ over 
$(\ZZ/p)[t^{\pm1}]$.

\section{$N(q,p)$-representations}

In this section, we discuss another type of
 metacyclic representations 
and the twisted Alexander polynomial associated to 
these representations.
Let $q \geq 1$ and $p = 2n+1$ be an odd prime. 
Consider a metacyclic group,
$N(q,p)= \ZZ/2q \marubackslash \ZZ/p$ 
that is a semi-direct product of $\ZZ/2q$ and
$\ZZ/p$ defined by 
\begin{equation}
N(q,p) = \langle
s,a| s^{2q}= a^p = 1 , sas^{-1}=a^{-1} 
\rangle.
\end{equation}


Note that $N(1,p)=D_p$ and 
$N(2,p)$ is a binary dihedral group, denoted by $N_p$. 
Since $s^2$
generates the center of $N(q,p)$, we see that 
$N(q,p)/\langle s^2\rangle =D_p$ 
and hence
$|N(q,p)| = 2pq$.
For simplicity, we assume hereafter that $\gcd(q,p)= 1$.
Now it is known \cite{Ha-M}, \cite{Ha}
that the knot group $G(K)$ of a knot $K$ is mapped
onto $N(q,p)$ if and only if $G(K)$ is mapped onto $D_p$, 
namely, $\Delta_K (-1) \equiv 0$ (mod $p$). 
For a 2-bridge knot $K(r)$, if 
$\Delta_{K(r)} (-1) \equiv 0$ (mod $p$), then we may
assume without loss of generality that there is an epimorphism
$\widetilde{\rho}: G(K(r))
\longrightarrow N(q,p)$ for any $q \geq 1$ such that
\begin{equation}
\widetilde{\rho}(x)= s\ {\rm  and}\ 
\widetilde{\rho}(y) = sa.
\end{equation}

As before, we draw a diagram below consisting of various groups and
connecting homomorphisms.

$
\begin{array}{
cc
ccc
cc}
& & & &GL(2qp,\ZZ)& &\\
& & &\hspace*{5mm} \mbox{\Large $\nearrow$}\mbox{{\large $\tilde{\xi}$}}& & &\\
& &\hspace*{5mm} N(q,p) & 
 \hspace*{5mm} \overset{\mbox{\large $\tilde{\pi}$}}{\longrightarrow}& 
 GL(2n,\CC)& 
 \overset{\mbox{\large $\tilde{\gamma}$}}{\longrightarrow}&GL(2nm,\ZZ)\\
 & \mbox{\large $\tilde{\rho}$}\mbox{\Large $\nearrow$}& & &&& \\
 G(K)& \overset{\mbox{{\large $\rho_p$}}}{\longrightarrow}&N_p&
 \overset{\mbox{\large $\xi_p$}}{\longrightarrow}&SU(2,\CC)& 
 \overset{\mbox{\large $\gamma_p$}}{\longrightarrow}&
 GL(4n,\ZZ)\\
 & \mbox{\large $\rho$}\mbox{\Large $\searrow$}&&
 \hspace*{5mm} &&&\\
 &&D_p&\hspace*{5mm}\overset{\mbox{{\large $\pi_0$}}}\longrightarrow&GL(2n,\ZZ)&&\\
 &&&\hspace*{5mm}  \mbox{\Large $\searrow$}\mbox{\large$\pi$}&&&\\
 &&&&GL(p,\ZZ)&&
\end{array}
$\\


Here, $p=2n+1, \hat{\rho}=\rho\circ\pi,\rho_0=\rho\circ\pi_0,
\tilde\nu=\tilde\rho\circ\tilde\xi,\tilde\tau=\tilde\rho\circ\tilde\pi$,$\tau_p=\rho_p\circ\xi_p$ 
and $m$ is the degree of the minimal polynomial of 
$\zeta$ over $\QQ$.
\\

Using the irreducible representation $\pi_0$ of $D_p$ on 
$GL(2n,\ZZ)$, we can
define a representation of $N(q,p)$ on $GL(2n,\CC)$. 
In fact, we have

\begin{lemm}\label{lem:newSec8.1}
Let $\zeta$ be a primitive $2q$-th root of $1$, $q \geq 1$. 
Then the mapping 
$\widetilde{\pi}: N(q,p) 
\longrightarrow GL(2n,\CC)$ defined by
\begin{align}
\widetilde{\pi} (s) &= \zeta \pi_0 (x)\ {\it and}\nonumber\\ 
\widetilde{\pi} (sa) &= \zeta \pi_0 (y)
\end{align}
gives a representation of $N(q,p)$ on $GL(2n,\CC)$.
\end{lemm}

Since a proof is straightforward, we omit details.
Now $\widetilde{\tau}= \widetilde{\rho} \circ \widetilde{\pi}:
G(K(r)) \longrightarrow GL(2n,\CC)$ defines a metacyclic
representation of $G(K(r))$. Then the twisted Alexander polynomial
$\widetilde{\Delta}_{\tilde{\tau}, K(r)}(t|\zeta)$
of $K(r)$ associated to $\widetilde{\tau}$ is given by
\begin{equation}
\widetilde{\Delta}_{\tilde\tau, K(r)}(t|\zeta) =
\widetilde{\Delta}_{\rho_0, K(r)}(\zeta t),
\end{equation}
where ${\rho}_0 =\rho \circ {\pi}_0$.

Therefore, the total twisted Alexander polynomial is
\begin{equation}
D_{\tilde{\tau}, K(r)}(t)= \prod_{(2q,k)=1}
\widetilde{\Delta}_{\rho_0, K(r)}({\zeta}^k t).
\end{equation}

This proves the following theorem.

\begin{thm}
Let $p=2n+1$ be an odd prime and $q \geq 1$.
Let $K(r)$ be a 2-bridge knot. 
Suppose $\Delta_{K(r)}(-1)\equiv 0$ (mod $p$).
Then $G(K(r))$ has a metacyclic representation
\begin{equation*}
\widetilde{\tau}=\widetilde{\rho} \circ \widetilde{\pi}: 
G(K(r)) \longrightarrow N(q,p) \longrightarrow GL(2n,\CC).
\end{equation*}

Let $\zeta$ be a primitive $2q$-th root of $1$. 
Then the twisted Alexander
polynomial $\widetilde{\Delta}_{\tilde{\tau}, K(r)}(t)$ and the total
twisted Alexander polynomial $D_{\tilde{\tau}, K(r)}(t)$
associated to $\widetilde{\tau}$ are given by
\begin{align}
&(1)\ \widetilde{\Delta}_{\tilde{\tau}, K(r)}(t)
=\widetilde{\Delta}_{\rho_0, K(r)}(\zeta t).\nonumber\\
&(2)\ D_{\tilde{\tau}, K(r)}(t)=
\prod_{(2q,k)=1}\widetilde{\Delta}_{\rho_0,K(r)}({\zeta}^k t).
\end{align}
\end{thm}
We conclude this section with a few remarks.
First, as we mentioned earlier, if $q=2$, $N(2,p)$ is a binary dihedral
group, denoted by $N_p$. 
It is known \cite{Kl} \cite{L2} 
that generators $s$ and $sa$ of $N_p$ are represented
in $SU(2, \CC)$ by trace free matrices. 
In fact, the mapping $\xi_p$:
\begin{equation}
\xi_p (s) =\mtx{0}{1}{-1}{0}\ 
{\rm and}\ 
\xi_p (sa) =\mtx{0}{v_p}{-v_p^{-1}}{0}
\end{equation}
gives a representation of $N_p$ into $SU(2, \CC)$, 
where $v_p=e^{\frac{2\pi i}{p}}$.

Then we will show that the total twisted Alexander polynomial
$D_{{\tau}_p, K(r)}(t)$
associated to ${\tau}_p = {\rho}_p \circ {\xi}_p$ is given by
\begin{equation}
D_{{\tau}_p, K(r)}(t)= \widetilde{\Delta}_{\rho_0,K(r)}(it) 
\widetilde{\Delta}_{\rho_0, K(r)}(-it),
\end{equation}
where $i = \sqrt{-1}$. 
Therefore we have the following corollary.
\begin{cor}
If $q=2$, then $D_{\tilde\tau, K(r)}(t) = 
D_{\tau_p, K(r)}(t)$.
\end{cor}
{\it Proof of (8.8).} 
Let $C_p$ be the companion matrix of the minimal
polynomial of $v_p$, namely,
$C_p=
\left[\begin{array}{ccc|c}
0&\cdots&0&-1\\\hline
& & &-1\\
& \hsymb{E}& &\vdots\\
& & &-1
\end{array}
\right].
$
Then, by definition, we have
\begin{equation}
D_{\tau_p, K(r)} (t) = \det[\widetilde{\Delta}_{\tau_p,K(r)}(t|C_p)].
\end{equation}
And (8.8) follows from the following lemma.
\begin{lemm}\label{lem:newSec8.4}
Let $E^*_{2n} = \left[a_{j,k}\right]$ be a 
$2n \times 2n$ matrix such that
$a_{j,k}= 1$, if $k+j =2n+1$ and $0$, otherwise
($E_{2n}^{*}$ is the \lq mirror image\rq\ of $E_{2n}$.)
Denote\\
\begin{align*}
&A = \mtxc{0}{E_{2n}}{-E_{2n}}{0}, 
B = \mtxc{0}{C_p}{-C_p^{-1}}{ 0},\ {\it and}\\
&\widehat A = \mtxc{i E_{2n}^{*}}{0}{0}{ -iE_{2n}^{*}}, 
\widehat B =\mtxc{i\pi_0(y)}{0}{0}{-i\pi_0(y)}.
\end{align*}
Then there exists a matrix $M_{4n} \in GL(4n,\CC)$ such that
$M_{4n}A {M_{4n}}^{-1}=\widehat A$ and $M_{4n} B {M_{4n}}^{-1} =
\widehat B$.
\end{lemm}

{\it Proof.}
 A simple computation shows that
$M_{4n}=\frac{1}{\sqrt{2}}\mtx{E_{2n}}{-iE_{2n}^{*}}{E_{2n}}{iE_{2n}^{*}}$
is what we sought. \qed

Secondly, the metacyclic group $N(q,p)$ is also represented by $\widetilde{\xi}$
in 
$GL(2qp,\ZZ)$ via
{\it maximum} permutation representation
on the symmetric group $S_{2qp}$.
To be more precise, let 
\begin{align*}
S=\{1,s,s^2,\cdots,s^{2q-1},\ 
& a,sa,s^2a,\cdots,s^{2q-1}a,\ 
a^2,sa^2,s^2a^2,\cdots,s^{2q-1}a^2,\ \cdots,\\
& a^{p-1},sa^{p-1},s^2a^{p-1},\cdots,s^{2q-1}a^{p-1}\}
\end{align*}
be the ordered set of the elements of $N(q,p)$.
Then the right
multiplication by an element $g$ of $N(q,p)$
on $S$ induces a permutation associated to $g$, and by taking the
permutation matrix corresponding to this permutation,
we obtain the representation $\widetilde\xi$ 
of $N(q,p)$ on $GL(2qp,\ZZ)$. 

Then we have the following:
\begin{prop}\label{prop:8.5}
For any $q\ge 1$, the twisted Alexander polynomial 
$\widetilde{\Delta}_{\tilde\nu,K(r)}(t)$ of $K(r)$
associated to $\widetilde{\nu}=\widetilde{\rho} \circ \widetilde{\xi}$
is given by
\begin{equation}
\widetilde{\Delta}_{\tilde{\nu}, K(r)}(t)
=
\frac{\displaystyle{\prod_{k=0}^{2q-1}} \Delta_{K(r)}(\zeta^k t)}{1-t^{2q}}
\prod_{k=0}^{2q-1} \widetilde{\Delta}_{\rho_0, K(r)}(\zeta^k t),
\end{equation}
where $\zeta$ is a primitive $2q$-th root of $1$.
Therefore,
$\widetilde{\Delta}_{\tilde{\nu},K(r)}(t)$ is an integer polynomial
in $t^{2q}$ and $D_{\tilde{\tau},K(r)}(t)$ divides
$\widetilde{\Delta}_{\tilde{\nu},K(r)}(t)$.
\end{prop}

{\it Proof.}
By construction, $\widetilde\xi(s)=\rho(x) \otimes C$ and
$\widetilde\xi(sa)=\rho(y) \otimes C$, where $C$ is the transpose of the
companion matrix of $t^{2q}-1$ and $[a_{i,j}] \otimes C=[a_{i,j}C]$,
the tensor product of $[a_{i,j}]$ and $C$.

Therefore (8.10) follows immediately. \qed

If Conjecture A holds for $K(r)$, $\widetilde{\Delta}_{\tilde{\nu},
K(r)}(t)$ is of the form:
\begin{equation*}
\widetilde{\Delta}_{\tilde{\nu}, K(r)}(t) =
\frac{\prod_{k=0}^{2q-1} \Delta_{K(r)}({\zeta}^k t)}{1-t^{2q}}
f(t^{2q})^2,
\end{equation*}
for some integer polynomial $f(t^{2q})$ in $t^{2q}$.

If coefficients are taken from a finite field, then (8.10)
becomes much simpler. The following proposition is a metacyclic version of
(2.12). Since a proof is easy, we omit details.

\begin{prop}
Let $p$ be an odd prime. Suppose $\Delta_{K(r)}(-1)\equiv 0$ (mod $p$). Then we have
\begin{equation}
\widetilde\Delta_{\tilde\nu,K(r)}(t)\equiv
\left\{\prod_{k=0}^{2q-1} \Delta_{K(r)}(\zeta^k t)\right\}^p / (1-t^{2q})^p\ {\rm (mod}\
p\ {\rm )}.
\end{equation}
\end{prop}

\begin{rem}  
In \cite{cha},
Cha defined the twisted Alexander invariant of a fibred
knot $K$ using its Seifert fibred surface. Evidently, this invariant is
closely related to our twisted Alexander polynomial.
For example, as is described in \cite[Example]{cha}, if we
consider a regular dihedral representation $\rho$, then the invariant he
defined is essentially
the same as the twisted Alexander polynomial associated to a regular
dihedral representation $\widetilde{\nu}=\widetilde{\rho} \circ
\widetilde{\xi}: G(K) \rightarrow D_p = N(1,p) \rightarrow
GL(2p,\ZZ)$ we discussed
in this section. More precisely, let $A_{\rho,K}(t)$ be Cha's twisted
Alexander
invariant associated to $\rho$ and
$\widetilde{\Delta}_{\widetilde{\nu},K}(t)$ the twisted Alexander
polynomial of
a knot $K$ associated to $\widetilde{\nu}$. Then we have
\begin{equation}
(1-t^2) \widetilde{\Delta}_{\widetilde{\nu},K}(t)=
A_{\rho,K}(t^2).
\end{equation}

We should note that $\widetilde{\Delta}_{\widetilde{\nu},K}(t)$ is an
integer polynomial in $t^2$. 
(See Proposition \ref{prop:8.5}) Details will appear elsewhere.
\end{rem}

\section{Example}
The following examples illustrate our main theorem.

\begin{ex}
Dihedral representations 
$\tau: G(K(r))\longrightarrow D_p\longrightarrow GL(2,\CC)$.

(I) Let $p=3$ and $n=1$. Then $\theta_1(z) = z+3$ and 
$\omega = -3$.\\
(a) $r=1/3$. $D_{\tau,K(1/3)}(t)=\widetilde{\Delta}_{\rho_0,K(1/3)}(t) = 1-t^2$.\\
(b) $r= 1/9$. $D_{\tau,K(1/9)}(t)=
\widetilde{\Delta}_{\rho_0,K(1/9)}(t)=(1-t^2)(1-t^3+t^6)(1 +t^3 +t^6)$.\\
(c) $r=5/27$. $D_{\tau,K(5/27)}(t)=\widetilde{\Delta}_{\rho_0,K(r)}(t)
=(1-t^2) (1+t -t^2
+t^3 +t^4) (1 -t -t^2 -t^3 +t^4)$.

Note $\Delta_{K(5/27)}(t)=(1 -t +t^2)
(2 - 2t +t^2 -2t^3 +2t^4)$ and\\
$2 - 2t +t^2 -2t^3 +2t^4 \equiv - (1 -t -t^2 -t^3 +t^4)$ (mod $3$), and
\begin{align*}
D_{\tau,K(5/27)}(t) 
&\equiv \dfrac{\Delta_{K(r)}(t)}{1+t} 
\dfrac{\Delta_{K(r)}(-t)}{1-t}\\
&\equiv 
\dfrac{(1+t)^2 (1 -t -t^2-t^3 +t^4)}{1+t}
\dfrac{(1-t )^2 (1+t -t^2+t^3 +t^4)}{1-t} \\
&\equiv (1-t^2) (1-t -t^2 -t^3 +t^4)
(1+t -t^2 +t^3 +t^4) \ 
{\rm (mod}\ 3).
\end{align*}

(II) Let $p=5$ and $n=2$. Then $\theta_2 (z) = z^2 + 5z +5$.\\
(a) $r=1/5$. $D_{\tau,K(r)}(t)=(1-t^2)^2
\Delta_{K(1/5)}(t) \Delta_{K(1/5)}(-t)$.\\
(b) $r= 19/85$ $D_{\tau,K(r)}(t)
=D_{\tau,K(1/5)}(t)f(t) f(-t)$,
where $f(t) = 1- 3t - 2t^2 + 4t^3 - t^4 - 4t^6 -3t^7 + 7t^8 - 3t^9
-4t^{10} -t^{12} +4t^{13} -2t^{14} - 3t^{15}+ t^{16}$, and
$\Delta_{K(r)}(t) = \Delta_{K(1/5)}(t) g(t)$,
where $g(t) = 2 -2t +2t^2 - 2t^3 + t^4 - 2t^5 
+ 2t^6 - 2t^7 + 2t^8$, and
$f(t) \equiv g(t)^2$ (mod $5$).

Since $\Delta_{K(1/5)}(t) \equiv (1+t)^4$ (mod $5$), 
we see
\begin{align*}
D_{\tau,K(r)}(t) 
&= D_{\tau,K(1/5)}(t) f(t) f(-t)\\
&=\left\{(1+t)^2 \Delta_{K(1/5)}(t) f(t) \right\} 
\left\{ (1-t)^2
\Delta_{K(1/5)}(-t) f(-t) \right\}\\
&\equiv \{(1+t)^6 g(t)^2\}\{(1-t)^6 g(-t)^2\}\\
&\equiv \{(1+t)^3 g(t)\}^2 \{(1-t)^3 g(-t)\}^2\\
&\equiv\left\{\dfrac{\Delta_{K(1/5)}(t) g(t)}{1+t}\right\}^2
\left\{\dfrac{\Delta_{K(1/5)}(-t) g(-t)}{1-t}\right\}^2\\
&\equiv \left\{\dfrac{\Delta_{K(r)}(t)}{1+t}\right\}^2
\left\{\dfrac{\Delta_{K(r)}(-t) }{1-t}\right\}^2\ 
{\rm (mod}\ 5).
\end{align*}

(c) $r=21/115$. $D_{\tau,K(r)}(t)= D_{\tau,K(1/5)}(t)f(t) f(-t)$,
where $f(t) = 4 +2t - 3t^2 - t^3 - 8t^5 - 3t^6+4t^7 + t^9 + 9t^{10}
+t^{11} +4t^{13} -3t^{14} -8t^{15} - t^{17} -3 t^{18} +2t^{19}
+4t^{20}$, and
$ \Delta_{K(r)}(t) = \Delta_{K(1/5)}(t) g(t)$,
where $g(t) = 2 -2t +2t^2 - 2t^3 +2t^4 -3t^5 +2 t^6 - 2t^7 + 2t^8 - 2t^9 +
2t^{10}$, and
$f(t) \equiv g(t)^2$ (mod $5$). Therefore, we see
\begin{equation*}
D_{\tau, K(r)}(t) \equiv \left\{
\dfrac{\Delta_{K(r)}(t)}{1+t}
\right\}^2
\left\{\dfrac{\Delta_{K(r)}(-t)}{1-t}\right\}^2\
{\rm (mod}\ 5).
\end{equation*}

\end{ex}

\begin{ex}\label{ex:9.2n}
Binary dihedral representations.  
$\tau_p:G(K(r))\longrightarrow N_p
\longrightarrow GL(2n,\CC)$

(I) Let $p=3$ and $n=1$.\\
(a)When  $r=1/9$, $D_{\tau_p, K(r)}(t)$ = $(1+t^2)^2(1-t^6+t^{12})^2$.\\
(b)When $r=5/27$, $D_{\tau_p, K(r)}(t) = (1+t^2)^2(1+3t^2+t^4+3t^6+t^8)^2$.

(II) Let $p=5$ and $n=2$.\\
(a) When $r=1/5$, $D_{\tau_p, K(r)}(t) = (1+t^2)^4(1-t^2+t^4-t^6+t^8)^2$.\\
(b) When $r =19/85$, $D_{\tau_p, K(r)}(t) = (1+t^2)^4(1-t^2+t^4-t^6+t^8)^2
f(t)^2$,
where $f(t)=
1+13t^2+26t^4+20t^6+13t^8+22t^{10}+40t^{12}
+33t^{14}+25t^{16}+33t^{18}+40t^{20}+22t^{22}
+13t^{24}+20t^{26}+26t^{28}+13t^{30}+t^{32}$.
\end{ex}

\begin{ex}\label{ex:9.3}
$N(q,p)$-representations.
$\widetilde{\nu}:G(K(r))
\longrightarrow N(q,p) 
\longrightarrow GL(2pq,\ZZ)$.

(I) Let $q=4,p=3,N(4,3)=\ZZ/8\marubackslash \ZZ/3$.\\
(a) $r=1/3$. $\widetilde\Delta_{\tilde{\nu},K(1/3)}(t)=
(1-t^8)(1+t^8+t^{16})$.\\
(b) $r=1/9$. 
$\widetilde\Delta_{\tilde{\nu},K(1/9)}=
(1-t^8)(1+t^8+t^{16})(1+t^{24}+t^{48})^3$.\\
(c) $r=5/27$.  
\begin{align*}
&\ \ \ \ \widetilde\Delta_{\tilde{\nu},K(5/27)}\\
&\ \ \ \ \ \ \ \  =
(1-t^8)(1+t^8+t^{16})(16+31t^{8}+16t^{16})^2(1-79t^8+129t^{16}
-79t^{24}+t^{32})^2.
\end{align*}

(II) Let $q=5,p=3, N(5,3)=\ZZ/10\marubackslash \ZZ/3$.\\
(a) $r=1/3$. $\widetilde\Delta_{\tilde{\nu},K(1/3)}(t)=
(1-t^{10})(1+t^{10}+t^{20})$.\\
(b) $r=1/9$. $\widetilde\Delta_{\tilde{\nu},K(1/9)}=
(1-t^{10})(1+t^{10}+t^{20})(1+t^{30}+t^{60})^3$.\\
(c) $r=5/27$.  
\begin{align*}
&\ \ \ \ \widetilde\Delta_{\tilde{\nu},K(5/27)}\\
&\ \ \ \ \ \ \ \ =
(1-t^{10})(1+t^{10}+t^{20})(1-228t^{10}-314t^{20}-228t^{30}+t^{40})^2\\
&\ \ \ \ \ \ \ \ \ \ \times
(1024+1201t^{20}+1024t^{40}).
\end{align*}

(III) Let $q=3,p=5$, $N(3,5)=\ZZ/6\marubackslash \ZZ/5$\\
(a) $r=1/5$. $\widetilde\Delta_{\tilde\nu,K(1/5)}(t)=
(1-t^6)^3(1+t^6+t^{12}+t^{18}+t^{24})^3$.\\
(b) $r=19/85$.
\begin{align*}
&\ \   \widetilde\Delta_{\tilde{\nu},K(19/85)}\\
&\ \ \ \   =
(1-t^6)^3(1+t^6+t^{12}+t^{18}+t^{24})^3\\
&\ \ \ \ \ \   \times 
(64+64t^6+48t^{12}+12t^{18}+49t^{24}+12t^{30}+48t^{36}+64t^{42}+64t^{48})\\
&\ \ \ \ \ \   \times 
(1-1243t^6+3335t^{12}+1570t^{18}-2423t^{24}+6320t^{30}-992t^{36}\\
&\ \ \ \ \ \ \ \  -2181t^{42} +9451t^{48}-2181t^{54}-992t^{60}+6320t^{66}-2423t^{72}\\
&\ \ \ \ \ \ \ \   +1570t^{78}+3335t^{84}-1243t^{90}+t^{96})^2
\end{align*}
\end{ex}

\section{$K$-metacyclic representations}

In this section, we briefly discuss $K$-metacyclic 
representations of the
knot group.
Let $p$ be an odd prime. Consider a group $G(p-1,p|k)$ that has the
following presentation:
\begin{equation}
G(p-1,p|k) = \langle s,a | s^{p-1}= a^p = 1, sas^{-1} =a^k
\rangle,
\end{equation}
where $k$ is a primitive $(p-1)$-st root of $1$ (mod $p$).

We call $G(p-1,p| k)$ a $K$-metacyclic group according to 
\cite{Fox}. 

\begin{prop}
Two $K$-metacyclic groups of the same order, $p(p-1)$ say,
are isomorphic.
\end{prop}

{\it Proof.} Let $G(p-1,p| \ell)=\langle
u,b | u^{p-1}= b^p = 1, u b u^{-1} =
b^{\ell}\rangle$ be
another $K$-metacyclic group. 
Since $\ell$ is also a primitive $(p-1)$-st
root (mod $p$),
we see that $\ell \equiv k^m$ (mod $p$), $1 \leq m \leq p-2$ for some $m$, 
where $m$ and $p-1$ are coprime.
Take two integers $\lambda$ and $\mu$ such that 
$m \lambda + (p-1) \mu= 1$. 
Then it is easy to show that a homomorphism 
$h:\ G(p-1,p|k) \rightarrow G(p-1,p|\ell)$ defined
by $h(s) = u^{\lambda}$ and $h(a) = b$ is in fact an isomorphism. \qed

The following proposition is also well-known.

\begin{prop} \cite{Fox}\cite{Ha}
Let $p$ be an odd prime. 
Suppose that $k$ is a primitive $(p-1)$-st root of 1 
(mod $p$). Then the knot group $G(K)$ is
mapped onto $G(p-1,p|k)$
if and only if $\Delta_K(k) \equiv 0$ (mod $p$).
\end{prop}

As is shown in \cite{Fox}, 
$G(p-1,p|k)$ is faithfully represented in $S_p$ by

\begin{equation}
\sigma (a) = (1 2 3 \cdots p)\  {\rm and}\ 
\sigma (s) = (k^{p-1} k^{p-2} \cdots k^2 k).
\end{equation}

Let $\pi_*:\  G(p-1,p |k) \rightarrow GL(p, \ZZ)$ 
be a matrix representation of
$G(p-1,p|k)$ via $\sigma$.
Now, let $K(r)$ be a $2$-bridge knot. 
Suppose that $\Delta_K(k) \equiv 0$
(mod $p$) for some
primitive $(p-1)$-st root of $1$ (mod $p$). 
Then a homomorphism $\delta:\  G(K(r)) \rightarrow
G(p-1,p | k)$ given by
\begin{equation}
\delta (x) = s\ {\rm  and}\ \delta (y)=sa,
\end{equation}
 induces a 
$K$-metacyclic representation 
$\Theta = \delta \circ \pi_* :\ 
G(K(r)) \rightarrow GL(p,\ZZ)$.

Then Conjecture A states that

\begin{equation}
\widetilde{\Delta}_{\Theta,K(r)} (t) =
\left[\frac{\Delta_{K(r)}(t)}{1-t}\right] F(t^{p-1}).
\end{equation}

We will see that (10,4) holds for the following knots 
including a non-$2$-bridge knot.

\begin{ex}
(1) Consider a trefoil knot $K$. 
Since $\Delta_K(-2) \equiv 0$ (mod 7) and 
$-2$ is a primitive $6$th root of 1 (mod 7),
$G(K)$ is mapped onto $G(6,7|-2)$. 
Then $(\delta \circ \sigma) (x)=
\sigma(s) = (132645)$ and
$(\delta \circ \sigma) (y) = \sigma (sa) = (146527)$ and we see
$\widetilde{\Delta}_{\Theta,K} (t)
= \left[\frac{\Delta_K(t)}{1-t}\right] (1-t^6)$.

(2) Let $K=K(1/9)$. Since $K(1/9) \in H(3)$, $G(K(1/9))$ is mapped onto
$G(6,7|-2)$, and

$\widetilde{\Delta}_{\Theta,K} (t) =\left[\frac{\Delta_K(t)}{1-t}\right]
(1-t^6)(1-t^6+t^{12})$.

(3) Let $K=K(5/27) \in H(3)$. Then

$\widetilde{\Delta}_{\Theta,K} (t) =\left[\frac{\Delta_K(t)}{1-t}
\right] (1-t^6)
(1-7t^6+9t^{12}-7t^{18}+t^{24})$.
\end{ex}

\begin{ex}
Consider a knot $K=K(5/9)$. 
Since $\Delta_K (t)=2 - 5t + 2t^2$,

$\Delta_K (2)$ = 0 and hence $G(K(5/9))$ is mapped onto $G(m,p|2)$ for
any odd prime $p$, where $m$ is a divisor of $p-1$.

If $p=5$ or $11$, 
then $2$ is a primitive $(p-1)$-st root of 1 (mod $p$). We
see then:

(i) For $p=5, \widetilde{\Delta}_{\Theta,K} (t) =
\left[\frac{\Delta_K(t)}{1-t}\right] (1-t^4)$.

(ii) For $p=11, \widetilde{\Delta}_{\Theta,K} (t) =
\left[\frac{\Delta_K(t)}{1-t}\right] (1-t^{10})$.

It is quite likely that we have $\widetilde{\Delta}_{\Theta,K} (t)
 =\left[\frac{\Delta_K(t)}{1-t}\right] 
 (1-t^{p-1})$, for any odd prime $p$ such that $2$
is a primitive $(p-1)$-st root of 1
(mod $p$).

(iii) If $p=7$,then $2$ is a primitive third root of $1$ 
(mod $7$) and hence
$G(K)$ has a representation $\Theta:\  
G(K) \rightarrow G(3,7|2)
\rightarrow GL(7,\ZZ)$ and we obtain

$\widetilde{\Delta}_{\Theta,K} (t) = 
\left[\frac{\Delta_K(t)}{1-t}\right] (1-t^3)^2$.
\end{ex}

\begin{ex}
Consider a non-2-bridge knot $K=8_5$ in 
Reidemeister-Rolfsen table. 
We have a Wirtinger presentation
$G(K)=\langle x,y,z|R_1,R_2\rangle$, where\\
\hspace*{2cm}
$R_1= 
(x^{-1}y^{-1}zyxy^{-1}x^{-1}y^{-1})
x(yxyx^{-1}y^{-1}z^{-1}yx)y^{-1}$ and\\
\hspace*{2cm}
$R_2 = (yx^{-1}y^{-1}z^{-1}x^{-1})y (xzyxy^{-1})z^{-1}$.
\hfill (10.5)
\setcounter{equation}{5}
Since $\Delta_{K} (t) = 
(1-t+t^2 )(1-2t+t^2 -2t^3 +t^4)$, it follows that
$\Delta_{K}(-1) \equiv 0$ (mod $3$) and 
$\Delta_{K} (-1) \equiv 0$ (mod $7$), and
further $\Delta_{K} (-2) \equiv 0$ (mod $7$). 
Therefore, $G(K)$ is mapped
onto each of the following groups: 
$D_3, D_7, N(2,3), N(2, 7)$ and $G(6,7|-2)$,
since $-2$ is a primitive $6$-th root of $1$ (mod 7).

Now we have five representations and computed their twisted Alexander
polynomials.

(1) For $\rho_1:\  G(K) \rightarrow D_3 \rightarrow GL(3,\ZZ)$, 
defined
by $\rho_1(x)= \rho_1(z)= \pi \rho (x)$ and 
$\rho_1(y) = \pi \rho(y)$,we have

$\widetilde{\Delta}_{\rho_1, K}(t) = 
\left[ \frac{\Delta_K(t)}{1-t}\right] f_1(t)
f_1(-t)$, where $f_1 (t)= (1+t)(1+t - 2t^2 +t^3 +t^4 )$.

(2) For $\rho_2:\ 
G(K) \rightarrow D_7 \rightarrow GL(7,\ZZ)$, defined by
$\rho_2(x)=\rho_2 (y)=\pi \rho (x)$ and 
$\rho_2(z) =\pi \rho(y)$, we have

$\widetilde{\Delta}_{\rho_2, K}(t) = 
\left[ \frac{\Delta_K(t)}{1-t}\right] f_2 (t)
f_2 (-t)$, where $f_2 (t)=(1+t)^3 (1+2t - 7t^3 -13t^4
-13t^5 -11t^6 -13t^7 -13t^8 -7t^9 +2t^{11} +t^{12})$.

(3) For $\rho_3:\ 
G(K) \rightarrow N(2,3) \rightarrow GL(12,\ZZ)$, defined
by $\rho_3(x)= \rho_3(z)= \widetilde{\nu} (x)$ and $\rho_3 (y)=
\widetilde{\nu} (y)$, we have
$\widetilde{\Delta}_{\rho_3, K}(t) = 
(1+t^2 )^2 (1+5t^2 +4t^4 +5t^6+t^8)^2$.

(4) For $\rho_4:\ G(K) \rightarrow N(2,7) \rightarrow GL(28,\ZZ)$, 
defined by $\rho_4(x)=\rho_4 (y)=\widetilde{\nu} (x)$ and $\rho_4
(z)=\widetilde{\nu} (y)$, we have
$\widetilde{\Delta}_{\rho_4, K}(t)=(1+t^2)^6 (1+4t^2 + 2t^4 +19t^6
+13t^8 +37t^{10} +17t^{12} +37t^{14} +13t^{16} +19t^{18} +2t^{20}
+4t^{22} +t^{24})^2$.

(5) For $\rho_5:\  
G(K) \rightarrow G(6,7|-2) \rightarrow GL(7,\ZZ)$, 
defined by $\rho_5(x)= \rho_5 (z)=\Theta (x)$ and 
$\rho_5 (y)=\Theta (y)$, we have
$\widetilde{\Delta}_{\rho_5, K}(t) = 
\left[ \frac{\Delta_K(t)}{1-t}\right]F(t)$,
where $F(t)=(1-t^6)(1-72t^6 -82t^{12} -72t^{18} +t^{24})$.
\end{ex}

We note that this example also supports Conjecture A.

\section{Appendix}
\noindent
{\bf 11.1. Proof of Proposition 2.1.}\\

Let $\theta_n (z) = c_0^{(n)} + c_1^{(n)} z + \cdots + c_n^{(n)}z^n$ be
the polynomial
defined in Section 2. 
Here $c_k^{(n)} =\binom{n+k}{2k} +2\binom{n+k}{2k+1}$.
Now we define four $n \times n$ 
integer matrices $A, A^*, B, B^*$ as follows:
$A=[A_{i,j}]$, where $A_{i,j} = a_{i,n-j+1}$,
$A^* = [A^{*}_{i,j}]$, where $A_{i,j}^{*} = - a_{i,j-1}$,
$B = [B_{i,j}]$, where $B_{i,j} = b_{i,n-j+1}$, and
$B^* = [B_{i,j}^*]$, where $B_{i,j}^* = b_{i,j}$.

Here $a_{j,k}$ and $b_{j,k}$ are given as follows.
\begin{align} 
&(1)\ a_{j,j} = b_{j,j} = 1\ {\rm  for}\ 1 \leq j \leq n.\nonumber\\
&(2)\ {\rm  For}\ 1 \leq j \leq k,
a_{j,k}=\binom{j+k-1}{2j-1}\ {\rm and}\ 
b_{j,k} = \binom{j+k-2}{2j-2}.\nonumber\\
&(3)\ {\rm If}\ 0 \leq k < j, a_{j,k} = b_{j,k} = 0.
\end{align}

\begin{lemm}\label{lemm:a1}
{\it The following formulas hold.}
\begin{align}
{\it For}\ &0 \leq k \leq n, \nonumber\\
&(1)\ c_k^{(n)} = a_{k+1,n+1} +
a_{k+1,n}.\nonumber\\
{\it For}\ &1 \leq j \leq k,\nonumber\\
&(2)\ b_{j,k} = a_{j,k} - a_{j,k-1},\nonumber\\
&(3)\ b_{j,k} = a_{j-1,k-1} + b_{j,k-1}\ {\it and}\nonumber\\
&(4)\ -2\sum_{k=j}^n b_{j,k} = a_{j-1,n} - c_{j-1}^{(n)} + b_{j,n}.
\end{align}
\end{lemm}

{\it Proof.} Only (4) needs a proof. 
Since $\sum_{k=j}^n b_{j,k} =
\sum_{k=j}^n (a_{j,k} - a_{j,k-1}) = a_{j,n}$,
we need to show that $-2a_{j,n} = a_{j-1,n} - c_{j-1}^{(n)} + b_{j,n}$.
However, it follows easily from (11.2) (1)-(3). 
\qed

Now these formulas are sufficient to show that the 
$2n \times 2n$ matrix $U_n =\mtx{A}{A^*}{B}{B^*}$
is what we sought. Since a proof is straightforward,
we omit the details.

\arraycolsep=1.7pt


\begin{ex}\label{ex:a2}
For $n=4, 5$, $U_n$ are given by
\begin{center}
$U_4=\left[\begin{array}{rrrrrrrr}
4& 3& 2& 1& 0& -1& -2& -3\\
10& 4& 1& 0& 0& 0& -1& -4\\
6& 1& 0& 0& 0& 0& 0& -1\\
1& 0& 0& 0& 0& 0& 0& 0\\
1& 1& 1& 1& 1& 1& 1& 1\\
6& 3& 1& 0& 0& 1& 3& 6\\
5& 1& 0& 0& 0& 0& 1& 5\\
1& 0& 0& 0& 0& 0& 0& 1
\end{array}\right]$, and
$U_5=\left[\begin{array}{rrrrrrrrrr}
5& 4& 3& 2& 1& 0& -1& -2& -3& -4\\
20& 10& 4& 1& 0& 0& 0& -1& -4& -10\\
21& 6& 1& 0& 0& 0& 0& 0& -1& -6\\
8& 1& 0& 0& 0& 0& 0& 0& 0& -1\\
1& 0& 0& 0& 0& 0& 0& 0& 0& 0\\
1& 1& 1& 1& 1& 1& 1& 1& 1& 1\\
10& 6& 3& 1& 0& 0& 1& 3& 6& 10\\
15& 5& 1& 0& 0& 0& 0& 1& 5& 15\\
7& 1& 0& 0& 0& 0& 0& 0& 1& 7\\
1& 0& 0& 0& 0& 0& 0& 0& 0& 1
\end{array}\right].
$
\end{center}
\end{ex}

\arraycolsep\labelsep

\noindent
{\bf
11.2. Proof of Lemma 5.2.
}

First we write down a solution $X=V_n$ of the equation 
$X^2 = 4E_n+C_n$.
Let us begin with the alternating Catalan series
\begin{equation}
\mu (y) = \sum_{k=0}^\infty b_k y^k,\ {\rm where}\ 
b_k=\frac{(-1)^k}{k+2} \binom{2k+2}{k+1}.
\end{equation}

Therefore, 
$\mu (y) = 1- 2y +5y^2 - 14y^3 + 132y^4 - 429y^5+ 1430y^6
+ \cdots$.
Let $\theta_n (z)=c_0^{(n)}+ c_1^{(n)} z+ \cdots + c_n^{(n)}z^n$
be the polynomial defined in
Section 2. 
Using $\theta_n (z)$, we define a new polynomial 
$f_n (x)=x^n \theta (x^{-1})=
a_0^{(n)} + a_1^{(n)} x + a_2^{(n)}x^2 + \cdots + a_n^{(n)} x^n$. 
For example, $f_1(x) = x \theta_1(x^{-1})=
x(3+x^{-1}) = 3x + 1$, and $f_2(x)= 5x^2+ 5x +1$.
Since $a_k^{(n)} = c_{n-k}^{(n)}$, we see that
\begin{equation}
a_k^{(n)}=
\frac{2n+1}{2n-2k+1} \binom{2n-k}{2n-2k}
=\binom{2n-k+1}{2n-2k+1} + \binom{2n-k}{2n-2k+1}.
\end{equation}

Next, we compute $f_n(x) \mu (y)$ = $\sum_{r,s \geq 0} c_{r,s}^{(n)} x^r
y^s$, where $c_{r,s}^{(n)} = a_r^{(n)}b_s$,
and define integers $d_{k,\ell}^{(n)}$, $0 \leq k, \ell$, as follows:
\begin{equation}
d_{k,\ell}^{(n)}=c_{k,\ell}^{(n)} + c_{k-1,\ell+1}^{(n)}
+ c_{k-2,\ell+2}^{(n)} + \cdots+ c_{0,k+ \ell}^{(n)}
={\displaystyle \sum_{i=0,i+j=k+\ell}^k} a_i^{(n)}b_j.
\end{equation}

Then we claim:

\begin{prop}\label{prop:a3}
$V_n$ = $[ v_{j,k}^{(n)}]_{1 \leq j,k \leq n}$, where
$v_{j,k}^{(n)} = d_{n-j,k-1}^{(n)}$, is a solution.
\end{prop}

\begin{ex}\label{ex:a4}
The following is the list of solutions $V_n, n=1,\dots,5$. \\
\arraycolsep=1.7pt
\begin{center}
$[1]$, 
$\mtx{3}{-5}{1}{-2}$,
$
\left[\begin{array}{rrr}
5&-7&14\\
5&-9&21\\
1&-2&5
\end{array}\right]$,
$
\left[\begin{array}{rrrr}
7& -9&18& -45\\
14&-23& 51&-132\\
7&-13& 31& -84\\
1& -2& 5& -14
\end{array}\right]$,
$
\left[\begin{array}{rrrrr}
9& -11& 22& -55& 154\\
30& -46& 99& -253& 715\\
27&-47&108&-286&825\\
9&-17& 41& -112& 330\\
1& -2& 5& -14& 42
\end{array}\right]$.
\end{center}
\end{ex}

\arraycolsep\labelsep

Now, to prove Proposition 
\ref{prop:a3}, we need several technical lemmas.
\begin{lemm}\label{lem:a5}
For $n \geq 2$ and $0 \leq k \leq n$, the following recursion
formula holds.
\begin{equation}
a_k^{(n)}=a_k^{(n-1)} + 2a_{k-1}^{(n-1)} -a_{k-2}^{(n-2)}.
\end{equation}
\end{lemm}

For convenience, we define $a_0^{(0)}$ = 1.
Since a direct computation using (11.4) verifies (11.6) easily, we omit
details.

Next, for $n, m \geq 0$, we define a number $F(n, m)$ as follows.
\begin{equation}
F(n, m) = \sum_{j=0}^{n}a_{n-j}^{(n)} b_{m+j}.
\end{equation}

\begin{ex}\label{ex:a6}
We have the following values for $F(n,m)$;
\begin{align*}
(1)\ &(i)\ F(0, 0) = a_0^{(0)} b_0= 1.\\
    &(ii)\ F(0, m) = a_0^{(n)}b_m = b_m.\\
(2)\ &(i)\ F(1,0) = a_1^{(1)}b_0 + a_0^{(1)} b_1 = 3 - 2 = 1.\\
&(ii)\ F(1,1) = a_1^{(1)}b_1 +a_0^{(1)}b_2 = -6 + 5 = -1.\\
&(iii)\ F(1, m) = a_1^{(1)} b_m + a_0^{(1)}b_{m+1} = 3 b_m +
b_{m+1}.\\
(3)\ &(i)\ 
F(2, 0) = a_2^{(2)}b_0 +a_1^{(2)}b_1 + a_0^{(2)}b_2 = 0.\\
&(ii)\ F(2, 1) = 1.\\
&(iii)\ F(2, 2) = - 3.
\end{align*}
\end{ex}

\begin{lemm}\label{lem:a7}
For $n \geq 2$ and $m \geq 0$, the following recursion formula
holds.
\begin{equation}
F(n, m) = F(n-1, m+1) + 2F(n-1, m) - F(n-2, m).
\end{equation}
\end{lemm}

\noindent
{\it Proof.} Use (11.6) to show (11.8) as follows:
\begin{align*}
F(n, m)&=\sum_{j=0}^{n}a_{n-j}^{(n)}b_{m+j}
        =\sum_{j=0}^{n}[a_{n-j}^{(n-1)} 
        + 2a_{n-1-j}^{(n-1)} -a_{n-2-j}^{(n-2)}]b_{m+j}\\
&=\sum_{j=0}^{n-1}a_{n-1-j}^{(n-1)}b_{m+1+j} +
2\sum_{j=0}^{n-1}a_{n-1-j}^{(n-1)} b_{m+j} -
\sum_{j=0}^{n-2}a_{n-2-j}^{(n-2)} b_{m+j}\\
&=F(n-1, m+1) + 2F(n-1, m) - F(n-2, m).\ \ \ \ \qed
\end{align*}

\begin{lemm}\label{lem:a8}
The following formulas hold.
\begin{align}
&(1)\ {\it For}\ n \geq 1\ {\it and}\ 
0 \leq k \leq n,
\sum_{j=0}^{n}a_{k-j}^{(k)}b_{j} = a_k^{(n-1)}.\nonumber\\
&(2)\ {\it For}\ n \geq 2\ {\it and}\ 
0 \leq m \leq n-2, F(n, m) =0.\nonumber\\
&(3)\ {\it For}\ n \geq 1, F(n, n- 1) = 1.\nonumber\\
&(4)\ {\it For}\ n \geq 1, F(n, n) = - (2n-1).
\end{align}
\end{lemm}

\noindent
{\it Proof.}
 (1) Use induction on $n$. Since (1) holds for $n = 1$, we may
assume
that it holds for $n$. Further, if $k = 0$, (1) holds trivially, and
hence it suffices
to show (1) for $n = n+1$ and $k = k+1$. Then, by (11.6),
\begin{align*}
\sum_{j=0}^{k+1}a_{k+1-j}^{(n+1)}b_j
&=\sum_{j=0}^{k+1}\{a_{k+1-j}^{(n)}
+ 2a_{k-j}^{(n)} -a_{k-1-j}^{(n-1)}\} b_j\\
&=\sum_{j=0}^{k+1}a_{k+1-j}^{(n)}b_j 
+ 2\sum_{j=0}^{n}a_{k-j}^{(n)} b_j
- \sum_{j=0}^{k-1}a_{k-1-j}^{(n-1)} b_j\\
&=a_{k+1}^{(n-1)} + 2a_k^{(n-1)} - a_{k-1}^{(n-2)}\\
& = a_{k+1}^{(n)}.
\end{align*}

\noindent
{\it Proof of (2).} 
Since 
 $F(n,m+1)=F(n+1,m)-2F(n,m)+F(n-1,m)$,
it suffices to show that $F(n,0)=0$ if $n\ge 2$.

Now
\begin{align*}
F(n, 0)&=
\sum_{j=0}^{n}a_{n-j}^{(n)} b_{j}
=
\sum_{j=0}^{n}\frac{(2n+1)}{(2j+1)}\binom{n+j}{2j}
\frac{(-1)^{j}}{(j+2)} \binom{2j+2}{j+1}\\
&= (2n+1) \sum_{j=0}^{n} (-1)^j \frac{(n+j)!}{(2j+1)!(n-j)!}
\frac{(2j+2)!}{(j+2)!(j+1)!}\\
&=(2n+1) \sum_{j=0}^{n} (-1)^j \frac{(n+j)! (2j+2)}{(n-j)!(j+2)!(j+1)!}\\
&= (2n+1) \sum_{j=0}^{n} (-1)^j\frac{2 (n+j)!}{(n-j)!(j+2)! j!}.
\end{align*} 
Therefore, to prove (2),
it suffices to show
\begin{equation}
\sum_{j=0}^{n} (-1)^j\frac{(n+j)!}{(n-j)!(j+2)! j!} =0
\end{equation}
or equivalently, by multiplying both sides through $n!/(n-2)!$, to show
\begin{equation}
\sum_{j=0}^{n} (-1)^j \binom{n}{j} 
\binom{n+j}{j+2} = 0.
\end{equation}

To show (11.11), we apply the following lemma
\cite[Lemma 5.3]{HM2}.


\begin{lemm}\label{lem:a9}
For $N \geq M \geq 0$ and $N \geq K \geq 0$,
\begin{align}
&\binom{N}{K} 
\binom{M}{M} - \binom{N-1}{K-1} \binom{M}{M-1} 
+\binom{N-2}{K-2} \binom{M}{M-2} - \cdots     \nonumber\\
&+(-1)^M \binom{N-M}{K-M} \binom{M}{0}         \nonumber\\
&= 
\binom{N-M}{K}.
\end{align}
\end{lemm}

Put $N = 2n$, $K = n+2$ and $M = n$ in (11.12). 
Since $N - M = n < K$, we see 
\begin{equation*}
\binom{2n}{n+2} \binom{n}{n}-\binom{2n-1}{n+1} 
\binom{n}{n-1}+ 
\cdots + (-1)^{n} \binom{n}{2} 
\binom{n}{0}
=\binom{n}{n+2}=0,
\end{equation*}
and hence $\sum_{j=0}^{n} (-1)^j 
\binom{n+j}{j+2} \binom{n}{j}=0$. 
This proves (2).

\medskip

{\it Proof of (3)}. By (11.8), we see that for $n \geq 2$,
\begin{equation*}
F(n+1, n-2) = F(n, n-1) + 2F(n, n-2) - F(n-1, n-2).
\end{equation*}
 Since $F(n+1, n-2)
=F(n, n-2) =0$ by (11.9) (2), it
follows that $F(n, n-1) = F(n-1, n-2)$ , and hence $F(n, n-1)=F(1, 0)= 1$
by Example \ref{ex:a6} (2)(i). 

{\it Proof of (4)}. Use (11.8) for $n \geq 1$ to see 
\begin{equation*}
F(n+1, n-1) = F(n, n) +
2F(n, n-1) - F(n-1, n-1).
\end{equation*}
Since $F(n+1, n-1)=0$ and $F(n, n-1)= 1$, it follows that
\begin{equation*}
F(n, n)=F(n-1, n-1) - 2
\end{equation*}
and hence,
\begin{equation*}
F(n, n)=F(1,1)-2(n-1)= -1-2n+2=-(2n-1).
\end{equation*}
\qed

We define another number $H_k^{(n)}$ as follows.
For any $n \geq 1$ and $k \geq 2$, we define
\begin{equation}
H_k^{(n)}=\sum_{j=0}^{k} a_j^{(n)} F(n-1,n+k-2-j) -
\sum_{j=0}^{k-2} a_j^{(n-1)} F(n,n+k-3-j).
\end{equation}
For example, $H_2^{(5)}=a_0^{(5)} F(4,5) + a_1^{(5)} F(4,4) +
a_2^{(5)} F(4,3) - a_0^{(4)} F(5,4)= 0$.\\
In particular, we should note; 
\begin{equation}
{\rm For}\ k \geq 2, H_k^{(1)} = 0.
\end{equation}

\begin{align*}
{\rm In\ fact},\ H_k^{(1)}&=a_0^{(1)} F(0, k-1) + a_1^{(1)} F(0, k-2)
-a_0^{(0)} F(1, k-2)\\
&=b_{k-1} + a_1^{(1)} b_{k-2} - 3b_{k-2} -b_{k-1}\\
&= 0.
\end{align*}

The last formula we need is the following lemma.

\begin{lemm}\label{lem:a10}
For any $n \geq 1$ and $k \geq 2$, we have
$H_k^{(n)} = 0.$ \hfill \mbox{{\rm (11.15)}}
\end{lemm}
\setcounter{equation}{15}

\noindent
{\it Proof.} 
We compute $H=H_k^{(n)} - H_k^{(n-1)}$. 
By definition, for $n \geq 2$,
\begin{align*}
H &=- \sum_{j=0}^{k-2} a_j^{(n-1)} F(n, n+k-3-j) + a_0^{(n)}F(n-1,
n+k-2)\\
&\ \  \ +a_1^{(n)} F(n-1, n+k-3)  
+ \sum_{j=0}^{k-2}(a_{j+2}^{(n)} +a_j^{(n-2)})F(n-1, n+k-4-j) \\
&\ \ \ -\sum_{j=0}^{k} a_j^{(n-1)} F(n-2, n+k-3-j).
\end{align*}

Since $a_{j+2}^{(n)} + a_j^{(n-2)}=a_{j+2}^{(n-1)} + 2a_{j+1}^{(n-1)}$
and $a_1^{(n)}=a_1^{(n-1)}+2a_0^{(n-1)}$ by (11.6),
we see
\begin{align*}
H  &=- \sum_{j=0}^{k-2} a_j^{(n-1)} F(n, n+k-3-j) + a_0^{(n)}F(n-1,
n+k-2) \\
&\ \ \ \ +(a_1^{(n-1)}+2a_0^{(n-1)})F(n-1, n+k-3)\\
&\ \ \ \ + \sum_{j=0}^{k-2}(a_{j+2}^{(n-1)} + 2a_{j+1}^{(n-1)} )F(n-1, n+k-4-j)\\
&\ \ \ \ -\sum_{j=0}^{k} a_j^{(n-1)} F(n-2, n+k-3-j).
\end{align*}

Note $a_0^{(n)} =a_0^{(n-1)}$ = 1 to see
\begin{align*}
H & =a_0^{(n-1)}\{-F(n, n+k-3) + F(n-1, n+k-2)\\
&\ \ \ \ \ \ \ \  +2F(n-1, n+k-3) - F(n-2,
n+k-3)\}\\
&\ \ \ \ + a_1^{(n-1)}\{-F(n, n+k-4) + F(n-1, n+k-3) \\
&\ \ \ \ \ \ \ \ +2F(n-1, n+k-4) - F(n-2,
n+k-4)\} +
\cdots \\
&\ \ \ \ 
+ a_{k-2}^{(n-1)}\{-F(n, n-1) + F(n-1, n) \\
&\ \ \ \ \ \ \ \ +2F(n-1, n-1) - F(n-2,
n-1)\}\\
&\ \ \ \ + a_{k-1}^{(n-1)}\{F(n-1, n-1) +2F(n-1, n-2) - F(n-2, n-2)\}\\
&\ \ \ \ 
+a_k^{(n-1)} \{ F(n-1, n-2)-F(n-2, n-3)\}.
\end{align*}

By (11.8) and (11.9)(3),(4), we see easily that each term of the summation
is
equal to $0$. This proves $H = 0$. 
\qed

Now we are in position to prove Proposition \ref{prop:a3}. 
Let $\mathbf{u}_j=(v_{j,1}^{(n)}, v_{j,2}^{(n)}, \cdots, v_{j,n}^{(n)})$
and $\mathbf{w}_k=(v_{1,k}^{(n)}, v_{2,k}^{(n)}, \cdots,
v_{n,k}^{(n)})^{T}$
be, respectively, the $j$-th row vector and the $k$-th column vector of $V_n$.
Then
we must show
\begin{align}
&(1)\ \mathbf{u}_{n-j} \cdot 
\mathbf{w}_k = 0\ {\it for}\ (i)\ 0 \leq j \leq
n-3, 1 \leq k \leq n-j-2\ {\it and}\ \nonumber\\
&\ \ \ \ 
(ii)\ 2\leq j \leq n-1, n-j+1 \leq k \leq n-1.\nonumber\\
&(2)\ 
\mathbf{u}_{n-j}\cdot \mathbf{w}_{n-j-1}= 1\ {\it  for}\ 
0 \leq j \leq n-2.\nonumber\\
&(3)\ 
\mathbf{u}_{n-j}\cdot \mathbf{w}_{n-j}= 4,\ 
{\it for}\ 1 \leq j \leq n-1.\nonumber\\
&(4)\ 
\mathbf{u}_n\cdot \mathbf{w}_n =
 4 - a_1^{(n)} = 4 - c_{n-1}^{(n)},\nonumber\\
&(5)\ 
\mathbf{u}_{n-j} \cdot 
\mathbf{w}_n = - a_{j+1}^{(n)} = - c_{n-j-1}^{(n)}\ 
{\it for}\ 1 \leq j \leq n-1.
\end{align}

Since (11.16) is obviously true for $n=1$, 
we assume hereafter that $n \geq2$.

We introduce new vectors, 
$\mathbf{b}_j=(b_j, b_{j+1}, \cdots,
b_{j+n-1})$ for $j \geq 0$ and $\mathbf{a}_k^{(n)}=
(a_k^{(n)}, a_{k-1}^{(n)}, \cdots, a_0^{(n)}, 0 , 
\cdots,0)^T$ for $0 \leq
k \leq n$. Then, from the definition
of $v_{j,k}^{(n)}$, it is easy to see the following:
\begin{align}
&(1)\ {\rm For}\ 
0 \leq j \leq n-1, 
\mathbf{u}_{n-j} = a_j^{(n)} \mathbf{b}_0
+ a_{j-1}^{(n)} \mathbf{b}_1 + \cdots + a_0^{(n)} 
\mathbf{b}_j.\nonumber\\
&(2)\ {\rm For}\ 
1 \leq k \leq n, \mathbf{w}_k = b_{k-1} 
\mathbf{a}_{n-1}^{(n)}
+ b_k \mathbf{a}_{n-2}^{(n)} + \cdots + b_{n+k-2} 
\mathbf{a}_0^{(n)}.
\end{align}

Since $\mathbf{u}_{n-j}\cdot \mathbf{w}_k=
\sum_{i=0}^{j} a_{j-i}^{(n)}
(\mathbf{b}_i \cdot \mathbf{w}_k)$, we first compute 
$\mathbf{b}_i \cdot
\mathbf{w}_k$.
In fact, a straightforward computation shows\\
\hspace*{1.5cm}$\mathbf{b}_i\cdot \mathbf{w}_k=
b_{k-1}( a_{n-1}^{(n)} b_i + 
a_{n-2}^{(n)} b_{i +1} + \cdots + a_0^{(n)} b_{n+i-1})\\
\hspace*{1.5cm}\hspace*{1.5cm}
+ b_k( a_{n-2}^{(n)} b_i + a_{n-3}^{(n)} b_{i +1} + \cdots + a_0^{(n)}b_{n+i-2})
+ \cdots + b_{n+k-2}(a_0^{(n)} b_i)$\\
\hspace*{1.5cm}\hspace*{1.1cm}
$= b_{k-1}(F(n, i-1) - a_n^{(n)} b_{i-1})\\
\hspace*{1.5cm}\hspace*{1.5cm}
+ b_k ( F(n, i-2) - a_{n-1}^{(n)} b_{i-1} - a_n^{(n)}
b_{i-2})\\
\hspace*{1.5cm}\hspace*{1.5cm}
+ \cdots\\
\hspace*{1.5cm}\hspace*{1.5cm}
+ b_{k+i-2} ( F(n, 0) - a_{n-i+1}^{(n)} b_{i-1}- \cdots -
a_n^{(n)} b_0)\\
\hspace*{1.5cm}\hspace*{1.5cm}
+b_{k+i-1}(a_{n-1}^{(n-1)} - a_{n-i}^{(n)} b_{i-1}- \cdots
- a_{n-1}^{(n)} b_0)\\
\hspace*{1.5cm}\hspace*{1.5cm}
+ \cdots\\
\hspace*{1.5cm}\hspace*{1.5cm}
 + b_{n+k-2}(a_i^{(n-1)} - a_1^{(n)} b_{i-1} - 
a_2^{(n)} b_{i-2} - \cdots - a_i^{(n)} b_0)\\
\hspace*{1.5cm}\hspace*{1.5cm}
+ b_{n+k-1}(a_{i-1}^{(n-1)} - a_0^{(n)} b_{i-1}- 
a_1^{(n)} b_{i-2} - \cdots - a_{i-1}^{(n)} b_0)\\
\hspace*{1.5cm}\hspace*{1.5cm}
+ \cdots \\
\hspace*{1.5cm}\hspace*{1.5cm}
+ b_{n+k+i-2}(a_0^{(n-1)} - a_0^{(n)} b_0)$.

Note that in the above summation, each of the last $i$ terms is
0 by (11.9)(1).

By rearranging this summation, we obtain\\
\hspace*{1.5cm}$\mathbf{b}_i \cdot \mathbf{w}_k=b_{k-1}F(n, i-1) + b_k F(n, i-2) + \cdots\\
\hspace*{1.5cm}\hspace*{1.5cm}
+ b_{k+i-2}F(n, 0) + F(n-1, k+i-1)\\
\hspace*{1.5cm}\hspace*{1.5cm} - b_{i-1}F(n, k-1) - b_{i-2}F(n, k) -
\cdots - b_0 F(n,k+i-2)$.

Since $0 \leq i \leq j \leq n-1$, we have for $\ell \geq 0$, $i-1- \ell
\leq n-2$ and hence $F(n, i-1- \ell)= 0$. Therefore
\begin{equation}
\mathbf{b}_i \cdot \mathbf{w}_k=F(n-1, k+i-1) -
\sum_{\ell=0}^{i-1} b_{i-1-\ell} F(n, k-1+\ell).
\end{equation}

Case 1. $i= 0$. Then $\mathbf{b}_0 \cdot \mathbf{w}_k=F(n-1, k-1)$. If 
$1\leq k \leq n-2$, then
$F(n-1, k-1)= 0$, and hence $\mathbf{u}_n \cdot \mathbf{w}_k=
a_0^{(n)}(\mathbf{b}_0 \cdot \mathbf{w}_k)= 0$.
Further, 
\begin{align*}
&\mathbf{u}_n \cdot \mathbf{w}_{n-1}=a_0^{(n)}F(n-1, n-2)=
a_0^{(n)} = 1,\ {\rm and}\\
&\mathbf{u}_n \cdot \mathbf{w}_n=a_0^{(n)}F(n-1, n-1)
= -(2n-3) =4-(2n+1) = 4 - a_1^{(n)}. 
\end{align*}
This proves (11.16) for $j= 0$.

Case 2. $i=1$. Then $\mathbf{b}_1 \cdot \mathbf{w}_k=F(n-1, k) - b_0
F(n,k-1)$.
If 1 $\leq k \leq n-3$, then $F(n-1, k) = F(n, k-1)= 0$ and 
$\mathbf{b}_1\cdot \mathbf{w}_k= 0$.
Since $\mathbf{b}_0 \cdot \mathbf{w}_k= 0$, we have 
$\mathbf{u}_{n-1} \cdot
\mathbf{w}_k= 0$ for $1 \leq k \leq n-3$.
Further, 
\begin{align*}
\mathbf{u}_{n-1} \cdot \mathbf{w}_{n-2}&=
a_1^{(n)}(\mathbf{b}_0 \cdot
\mathbf{w}_{n-2}) + a_0^{(n)} (\mathbf{b}_1 \cdot \mathbf{w}_{n-2})\\
&=a_1^{(n)}F(n-1, n-3) + a_0^{(n)} \{F(n-1, n-2) - b_0 F(n, n-3)\}\\
&=a_0^{(n)} F(n-1, n-2) = 1.
\end{align*}
Also,
\begin{align*}
\mathbf{u}_{n-1} \cdot \mathbf{w}_{n-1} &= a_1^{(n)} F(n-1, n-2) +
a_0^{(n)} \{F(n-1, n-1)- b_0 F(n, n-2)\}\\
&= a_1^{(n)} + a_0^{(n)}(-(2n-3))= 2n+1 -(2n-3) = 4.
\end{align*}
Finally,
\begin{align*}
\mathbf{u}_{n-1} \cdot \mathbf{w}_n&=a_1^{(n)}F(n-1, n-1) +
a_0^{(n)} \{F(n-1, n) - b_0 F(n, n-1)\}\\
&=H_2^{(n)} - a_2^{(n)}=- a_2^{(n)},\ {\rm  by\ (11.15)}. 
\end{align*}
This proves (11.16)
for $j=1$.

Now we assume that $2 \leq j \leq n-1$ and compute 
$\mathbf{u}_{n-j} \cdot
\mathbf{w}_k$, $1 \leq k \leq n$. Then
\begin{align*}
\mathbf{u}_{n-j} \cdot \mathbf{w}_k
&=\sum_{i=0}^{j} a_{j-i}^{(n)} (
\mathbf{b}_i \cdot \mathbf{w}_k)\\
&= a_j^{(n)} (\mathbf{b}_0 \cdot \mathbf{w}_k)
+ \sum_{i=1}^{j} a_{j-i}^{(n)} (\mathbf{b}_i \cdot \mathbf{w}_k)\\
&=a_j^{(n)}F(n-1, k-1) + \sum_{i=1}^{j}a_{j-i}^{(n)} \{F(n-1, k+i-1) \\
&\ \ \ \ -\sum_{\ell=0}^{i-1} b_{i-1-\ell}F(n, k-1+ \ell)\}\\
&= \sum_{i=0}^{j} a_{j-i}^{(n)} F(n-1, k+i-1) -
\sum_{i=1}^{j}a_{j-i}^{(n)} \sum_{\ell=0}^{i-1} b_{i-1-\ell} F(n,
k-1+\ell).
\end{align*}

Therefore, the coefficient of $F(n, k-1+q)$, $0 \leq q \leq i-1$, is equal to
\begin{align*}
\sum_{i=1}^{j} a_{j-i}^{(n)} b_{i-1-q}
&=\sum_{i=q+1}^{j}a_{j-i}^{(n)} b_{i-1-q}\\
&=a_{j-q-1}^{(n-1)}\ {\rm by}\ (11.9)(1),\ {\rm and\ hence}
\end{align*}
\begin{equation}
\mathbf{u}_{n-j} \cdot \mathbf{w}_k=\sum_{i=0}^{j} a_{j-i}^{(n)}
F(n-1, k+i-1) - \sum_{q=0}^{j-1}a_{j-q-1}^{(n-1)}F(n, k-1+q).
\end{equation}

If $1 \leq k \leq n-j-2$, then $k+i-1 \leq k+j-1 \leq n-3$ and also $k-1+q
\leq k+j-2 \leq n-4$,
and hence, $\mathbf{u}_{n-j} \cdot \mathbf{w}_k = 0$.

If $k = n-j-1$, then $\mathbf{u}_{n-j} \cdot \mathbf{w}_{n-j-1}=
a_0^{(n)}F(n-1, n-2)=a_0^{(n)}= 1$. 
Further,
$\mathbf{u}_{n-j} \cdot \mathbf{w}_{n-j}=a_0^{(n)} F(n-1, n-1) + a_1^{(n)}
F(n-1, n-2)=-(2n-3) + 2n+1= 4$.

Now suppose $n-j+1 \leq k \leq n-1$. Then by (11.9),
\begin{align*}
\mathbf{u}_{n-j} \cdot \mathbf{w}_k 
&=\sum_{i=0}^{j}a_{j-i}^{(n)} F(n-1,k+i-1) 
- \sum_{q=0}^{j-1}a_{j-q-1}^{(n-1)} F(n,k-1+q)\\
&= \sum_{i=n-k-1}^{j}a_{j-i}^{(n)} F(n-1, k+i-1) - 
\sum_{q=n-k}^{j-1}a_{j-q-1}^{(n-1)} F(n,k-1+q).
\end{align*}

That is exactly $H_{k-(n-j-1)}^{(n)}$ and hence $\mathbf{u}_{n-j} \cdot
\mathbf{w}_k= 0$ for $n-j+1 \leq k \leq n-1$.

Finally, a similar computation shows that
\begin{align*}
\mathbf{u}_{n-j} \cdot \mathbf{w}_n&=\sum_{i=0}^{j} a_{j-i}^{(n)} F(n-1,
n+i-1) - \sum_{q=0}^{j-1}a_{j-1-q}^{(n-1)} F(n,n+q-1)\\
&=H_{j+1}^{(n)} - a_{j+1}^{(n)} F(n-1,n-2)\\
&=- a_{j+1}^{(n)}\\
&= -c_{n-j-1}^{(n)}.
\end{align*}

This proves (11.16) and a proof of Proposition \ref{prop:a3} 
is now complete.
\qed

\medskip
\medskip\medskip
\noindent
{\bf Acknowledgements. } 
The first author is
partially supported by MEXT, Grant-in-Aid for
Young Scientists (B) 18740035,
and the second author is
partially supported by NSERC Grant~A~4034.
The authors thank Daniel Silver and Joe Repka
for their invaluable comments.

\end{document}